\documentclass[preprint, 12pt]{elsarticle}

\usepackage{hyperref, lineno}
\usepackage{amsmath, amssymb}
\usepackage{amsfonts}
\usepackage{graphicx}
\usepackage{bm, bbm}
\usepackage{algorithm, caption}
\usepackage{algpseudocode}
\usepackage{afterpage}
\usepackage{physics}
\usepackage{subcaption}
\usepackage{multirow}
\usepackage{color}
\usepackage{lineno}
\usepackage{mathtools}
\usepackage{hhline}

\newcommand{\Cov}{\text{Cov}}
\newcommand{\R}{\mathbb{R}}

\journal{Journal Name}

\begin{document}

\begin{frontmatter}

\title{
Quantifying total uncertainty in physics-informed neural networks for solving forward and inverse stochastic problems
}

\author[Brown]{Dongkun Zhang}
\author[Brown]{Lu Lu}
\author[SNH]{Ling Guo\corref{cor}} \ead{lguo@shnu.edu.cn}
\author[Brown]{George Em Karniadakis}
\address[Brown]{Division of Applied Mathematics, Brown University, Providence RI, USA}
\address[SNH]{Department of Mathematics, Shanghai Normal University, Shanghai, China}
\cortext[cor]{Corresponding Author}

\begin{abstract}
Physics-informed neural networks (PINNs) have recently emerged as an alternative way of solving partial differential equations (PDEs) without the need of building elaborate grids, instead, using a straightforward implementation. In particular, in addition to the deep neural network (DNN) for the solution, a second DNN is considered that represents the residual of the PDE. The residual is then combined with the mismatch in the given data of the solution in order to formulate the loss function. This framework is effective but is lacking uncertainty quantification of the solution due to the inherent randomness in the data or due to the approximation limitations of the DNN architecture. Here, we propose a new method with the objective of endowing the DNN with uncertainty quantification for both sources of uncertainty, i.e., the {\em parametric uncertainty} and the {\em approximation uncertainty}. We first account for the parametric uncertainty when the parameter in the differential equation is represented as a stochastic process. Multiple DNNs are designed to learn the modal functions of the {\em arbitrary polynomial chaos} (aPC) expansion of its solution by using stochastic data from sparse sensors. We can then make predictions from new sensor measurements very efficiently with the trained DNNs. Moreover, we employ {\em dropout} to correct the over-fitting and also to quantify the uncertainty of DNNs in approximating the modal functions. We then design an {\em active learning} strategy based on the dropout uncertainty to place new sensors in the domain in order to improve the predictions of DNNs. Several numerical tests are conducted for both the forward and the inverse problems to quantify the effectiveness of PINNs combined with uncertainty quantification. This NN-aPC new paradigm of physics-informed deep learning with uncertainty quantification can be readily applied to other types of stochastic PDEs in multi-dimensions.

\end{abstract}


\end{frontmatter}


\section{Introduction}
\label{S:1}
How to make the best use of existing data while exploiting the information from classical mathematical models or even empirical correlations developed within a discipline is an important issue as {\em data-driven modeling} is emerging as a powerful paradigm for physical and biological systems. For example, in geophysics, researchers have been using the remote sensing data collected from multi-spectral satellites and the top-of-atmospheric reflectance model as a calibration of the data to study the soil salinization~\cite{soil_sanilization}, or estimating the Earth heat loss based on the heat flow measurements and a model of the hydrothermal circulation in the oceanic crust~\cite{pollack1993heat}. Data can be used to provide closures in nonlinear models or to estimate parameters or functions in mathematical models. Moreover, mathematical models can be used as additional knowledge to formulate ``informative priors" in statistical estimation methods or be encoded in specially designed machine learning tools so that a smaller amount of data is required for inference of system identification. There has been recent progress for both forward (inference) and inverse (identification) problems using different methods. For example, for the forward problem, some of the popular choices of machine learning tools are Gaussian process~\cite{graepel2003, sarkka2011, bilionis2016, Raissi_nonlinear, Pang2018} and deep neural networks (DNNs)~\cite{Lagaris1997, Lagaris2000, Yinglexing17, MaziarParisGK17_1, Hadi18}. For inverse problems, similar methods have been advanced, e.g., Bayesian estimation~\cite{Stuart10} and variational Bayes inference~\cite{zhu_zabaras18}, and have been proposed for a wide variety of objectives, from parameter estimation~\cite{raissi_jcp_2017} to discovering partial differential equations~\cite{Rudy17, RudyandBrunton18, MaziarParisGK17_2, MaziarGK18JCP} to learning constitutive relationships~\cite{Tartakovsky2018}.

In this work we focus on the DNNs, and in particular the {\em physics-informed neural networks} (PINNs) for forward and inverse problems, first introduced in~\cite{MaziarParisGK17_1, MaziarParisGK17_2}. However, in those works the mathematical models were deterministic differential equations, so here we consider stochastic differential equations that model either random micro-structure in a medium or the lack of complete knowledge (``uncertainty"), e.g., of the material property. There have only been very few works published on solving stochastic differential equations using DNNs, e.g., \cite{Weinan-arxiv, Maziar-arxiv}, for forward problems. Here we study the special and perhaps the most complex case where some of the physics is known, namely via the stochastic differential equations, and the parameter in the equation is represented as a stochastic process, introducing {\em parametric uncertainty}. First, we solve the forward stochastic Poisson equation, where there is uncertainty associated with the driving force. Subsequently, we consider the inverse stochastic elliptic equation, where the diffusivity is modeled as a random process. In the latter case, we have only partial information of the diffusivity from scattered sensors but we have much more data available for the solution, and we aim to infer the stochastic processes of not only the solution but also the diffusivity, and quantify their uncertainties given the randomness in the data. An additional uncertainty is due to the DNN approximation, which we will refer to as the {\em approximation uncertainty}. Taken together, we refer to the parametric uncertainty and the approximation uncertainty as the {\em total} uncertainty. To the best of our knowledge, the current work is the first to address total uncertainty in solving stochastic forward and inverse problems using DNNs. 

In particular, in this paper we combine the {\em arbitrary polynomial chaos} (aPC) with PINNs for both the forward and the inverse stochastic problems. One of the most popular methods for uncertainty quantification studies is the {\em polynomial chaos}~\cite{Ghanem, XiuGK02} because it has been very effective in representing correlated stochastic fields. However, aPC~\cite{Summer-Wan-GK, WanGK06, Nowak12, Leihuan2018, Witteveen06} is more suitable for building the orthogonal basis from arbitrary random space, without the need of any assumption on the distribution of the data. Therefore, in the current work, we employ the aPC to develop a combined method that we call NN-aPC, where we use the DNNs to learn each individual mode of the aPC expansion. More importantly, after training, the proposed method can be used to predict new realizations of the solution based only on very few measurements.
 
Treatment of the DNN approximation uncertainty has been addressed using different methods in the past. The traditional way to estimate uncertainty in DNNs is using the Bayes' theorem, e.g., the Bayesian neural networks (BNNs)~\cite{mackay1992practical,neal2012bayesian}. BNNs are standard DNNs with prior probability distributions placed over their weights, and given observed data inference is then performed on weights. Because the inference is not tractable in general, variational inference is often used to approximate the inference~\cite{jordan1998introduction, paisley2012variational, kingma2013auto, rezende2014stochastic, titsias2014doubly, hoffman2013stochastic}. However, these models have very high additional computational cost because they require more parameters for the same network size and more time for the DNN parameters to converge. Recently, Gal et al. developed a new way to quantify uncertainty in DNNs by using {\em dropout}~\cite{Gal2016a, Gal2017a, Gal2017b}, which is largely used as a regularization technique~\cite{hinton2012improving, Srivastava2014} to address the problem of over-fitting. Gal et al.~\cite{Gal2016a} showed that a DNN with dropout is mathematically equivalent to approximating a probabilistic deep Gaussian process~\cite{damianou2013deep}, no matter what network architecture and non-linearities are used. Moreover, dropout does not induce much computation overhead and thus has been used as a practical tool to obtain uncertainty estimation effectively in real applications including language modelling~\cite{Gal2016b}, computer vision~\cite{kendall2015bayesian, Gal2017c} and medical applications~\cite{angermueller2017deepcpg, yang2016fast}. In this paper, dropout is used to to estimate the uncertainty in approximating each aPC mode. Based on the magnitude of this uncertainty we set up an active learning strategy and deploy additional sensors to obtain more measurements of the quantity of interest (QoI), in order to improve the predictability of PINNs.

The organization of this paper is as follows. In Section \ref{S:2}, we set up the data-driven forward and inverse problems. In Section \ref{S:3}, we introduce the PINNs for solving deterministic differential equations, followed by our main algorithm, the NN-aPC, and the method of dropout for uncertainty. In Section \ref{S:4}, we provide a detailed study of the accuracy and performance of the NN-aPC method for solving both the forward and inverse stochastic diffusion equation and demonstrate the effectiveness of active learning via dropout-induced uncertainty. Finally, we conclude with a brief discussion in Section \ref{S:5}.

\section{Problem Setup}
\label{S:2}
Suppose we have a stochastic differential equation:
\begin{equation}\label{eqn:SDE}
\begin{gathered}
    \mathcal{N}_x[u(x;\omega);k(x;\omega)] = 0, \quad x\in\mathcal{D}, \quad \omega\in\Omega,\\
    \text{B.C.:}\qquad \mathcal{B}_x[u(x;\omega)] = 0, \quad x\in\Gamma,
\end{gathered}
\end{equation}
where $\mathcal{N}_x$ is the general form of a differential operator that could be nonlinear, $\mathcal{D}$ is a $d$-dimensional physical domain in $\R^d$, $\Omega$ is the random space, and $u(x;\omega)$ is the solution to this equation. The boundary condition is imposed through the generalized boundary condition operator $\mathcal{B}_x$ at the domain boundary $\Gamma$. The random parameter $k(x;\omega)$ is the source of parametric uncertainty, which could be represented by either a few random variables or by an infinite dimensional (in the random space) random process. 

We consider two types of problems here: first, a {\em forward} problem, where we know exactly the distribution of $k(x;\omega)$ everywhere in the domain $\mathcal{D}$ and $u(x;\omega)$ is our QoI; and second, an {\em inverse} problem, where we assume that we have incomplete information on $k(x;\omega)$ but some extra knowledge on $u(x;\omega)$, and we are interested in inferring the full stochastic profile of $k(x;\omega)$. In practice, both problems are data-driven, since the information usually comes from data collected via sensor measurements. Here, we summarize the different scenarios of the sensors placement for each type of the problems:
\begin{itemize}
    \item \textit{Forward problem}: The $u$-sensors are placed only at the boundary $\Gamma$ to provide boundary condition, while the $k$-sensors are virtual (since we know the distribution of $k$), thus we can have as many $k$-sensors as we want and they can be placed anywhere in $\mathcal{D}$.
    \item \textit{Inverse problem}: In addition to having $u$-sensors at the boundary $\Gamma$, we have a limited number of extra $u$-sensors that can be placed in the domain $\mathcal{D}$, whereas we only have a limited number of $k$-sensors.
\end{itemize}

In this paper, we address both types of problems but we will focus more on solving the inverse problem.

\section{Methodology}
\label{S:3}
\subsection{Physics-informed neural network}
\label{S:3-1}
In this part, we briefly review using DNNs to solve deterministic differential equations~\cite{Lagaris1997, Lagaris2000, MaziarParisGK17_1}, and its generalization for solving deterministic inverse problems in~\cite{MaziarParisGK17_2}. To this end, re-consider Equation \ref{eqn:SDE} but replace the random input $\omega$ and approximate it with a finite set of parameters, leading to a parameterized differential equation:
\begin{equation}\label{eqn:ODE}
\begin{gathered}
    \mathcal{N}_x[u;\eta]=0,\quad x\in \mathcal{D},\\
    \text{B.C.:}\qquad \mathcal{B}_x[u] = 0, \quad x\in\Gamma,
\end{gathered}
\end{equation}
where $u(x)$ is the solution and $\eta$ denotes the parameters.

\begin{figure}[htbp]
	\centering
	\includegraphics[width=0.6\linewidth]{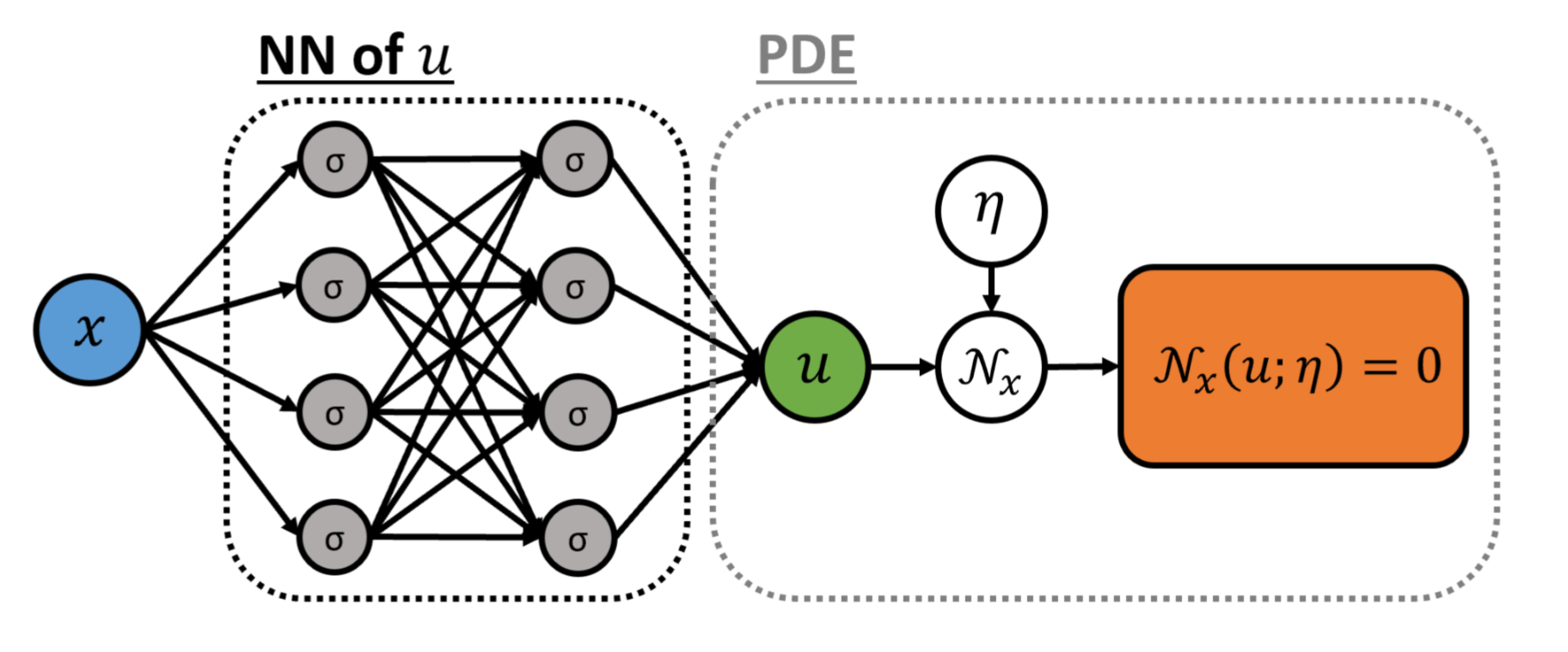}
	\caption{Schematic of the PINN for solving differential equations.}
	\label{fig:nn_sketch}
\end{figure}

A DNN, denoted by $\hat{u}(x;\theta)$, is constructed as a surrogate of the solution $u(x)$, and it takes the coordinate $x$ as the input and outputs a vector that has the same dimension as $u$. Here we use $\theta$ to denote the DNN parameters that will be tuned at the training stage, namely, $\theta$ contains all the weights $\bm{w}$ and biases $\bm{b}$ in $\hat{u}(x;\theta)$. For this surrogate network $\hat{u}$, we can take its derivatives with respect to its input by applying the chain rule for differentiating compositions of functions using the automatic differentiation, which is conveniently integrated in many machine learning packages such as Tensorflow~\cite{tensorflow}. The restrictions on $\hat{u}$ is two-fold: first, given the set of scattered data of the $u(x)$ observations, the network should be able to reproduce the observed value, when taking the associated $x$ as input; second, $\hat{u}$ should comply with the physics imposed by Equation \ref{eqn:ODE}. The second part is achieved by defining a residual network:
\begin{equation}\label{eqn:PINN}
    \hat{f}(x;\theta,\eta):= \mathcal{N}_x[\hat{u}(x;\theta);\eta],
\end{equation}
which is computed from $\hat{u}$ straightforwardly with automatic differentiation. This residual network $\hat{f}$, also named the physics-informed neural network (PINN), shares the same parameters $\theta$ with network $\hat{u}$ and should output the constant 0 for any input $x\in\mathcal{D}$. Figure \ref{fig:nn_sketch} shows a sketch of the PINN. At the training stage, the shared parameters $\theta$ (and also $\eta$, if it is also to be inferred) are fine-tuned to minimize a loss function that reflects the above two constraints.

Suppose we have a total number of $N_u$ observations on $u$, collected at location $\{x_u^{(i)}\}_{i=1}^{N_u}$, and $N_c$ is the number of collocation points $\{x_f^{(i)}\}_{i=1}^{N_c}$ where we evaluate the residual $\hat{f}(x_f^{(i)}; \theta, \eta)$. We shall use $(x^*, y^*)$ to represent a single instance of training data, where the first entry $x^*$ denotes the input and the second entry $y^*$ denotes the anticipated output (also called ``label''). The workflow of solving a differential equation with PINN can be summarized as follows:
\begin{algorithm}[H]
\caption{PINN for solving differential equations with random inputs}\label{alg:PINN}
\begin{algorithmic}[]
\State\textbf{Step 1:} Specify the training set:
$$\hat{u}\text{ network: }\{(x_u^{(i)}, u(x_u^{(i)}))\}_{i=1}^{N_u}, \quad \hat{f}\text{ network: }\{(x_f^{(i)}, 0)\}_{i=1}^{N_f}.$$
\State\textbf{Step 2:} Construct a DNN $\hat{u}(x;\theta)$ with random initialized parameters $\theta$.
\State\textbf{Step 3:} Construct the residual network $\hat{f}(x;\theta,\eta)$ by substituting the surrogate $\hat{u}$ into the governing equation (Equation \ref{eqn:PINN}) via automatic differentiation and arithmetic operations.
\State\textbf{Step 4:} Specify a loss function by summing the mean squared error of both the $u$ observations and the residual:
\begin{equation}\label{eqn:PINN_loss}
  \mathcal{L}(\theta,\eta)=\frac{1}{N_u}\sum_{i=1}^{N_u}[\hat{u}(x_u^{(i)};\theta)-u(x_u^{(i)})]^2+\frac{1}{N_c}\sum_{i=1}^{N_c}\hat{f}(x_f^{(i)};\theta,\eta)^2.
\end{equation}
\State\textbf{Step 5:} Train the DNNs to find the best parameters $\theta$ and $\eta$ by minimizing the loss function:
\begin{equation}\label{PINN_NN_training}
  \theta=\arg\min \mathcal{L}(\theta,\eta)
  \end{equation}
\end{algorithmic}
\end{algorithm}

\subsection{NN-aPC: Combining arbitrary polynomial chaos with neural networks}
\label{S:3-2}
We generalize the PINN method to solve stochastic differential equations for both forward and inverse problems, i.e., we aim to infer continuous random processes. Assume that a sensor will generate a sequence of measurements after being installed, and when the data is recorded, all sensors are read simultaneously. We denote the measurements from all the sensors at the same instant by a {\em snapshot} of the sensor data. Although the measurement results change from one measurement to the next due to randomness, it is reasonable to believe that every snapshot of sensor data corresponds to the same random event in the random space. We also assume that when the number of snapshots is big enough, the empirical distribution approximates the true distribution.

Let us consider Equation \ref{eqn:SDE}. Suppose we have $N_k$ sensors for $k(x;\omega)$ placed at $\{x_k^{(i)}\}_{i=1}^{N_k}$, $N_u$ sensors for $u(x;\omega)$ placed at $\{x_u^{(i)}\}_{i=1}^{N_u}$, and $N_f$ collocation points at $\{x_f^{(i)}\}_{i=1}^{N_f}$ that are used to calculate the residual of Equation \ref{eqn:SDE}. A total number of $N$ snapshots of measurements are made from all these sensors. Let $k_s^{(i)}$ and $u_s^{(i)}$ ($s=1,2,...,N$) be the $s$-th measurement of $k$ and $u$ at location $x_k^{(i)}$ and $x_u^{(i)}$ respectively, and $\omega_s$ is the random instance at the $s$-th measurement, i.e., $k_s^{(i)} = k(x_k^{(i)};\omega_s)$ and $u_s^{(i)} = u(x_u^{(i)};\omega_s)$. The training data set can be represented by
\begin{equation}\label{eqn:training_set}
    \mathcal{S}_t = \{\{(x_k^{(i)}, k_s^{(i)})\}_{i=1}^{N_k}, \{(x_u^{(i)}, u_s^{(i)})\}_{i=1}^{N_u}, \{(x_f^{(i)}, 0)\}_{i=1}^{N_f}\}_{s=1}^N.
\end{equation}
The proposed NN-aPC method consists of the following steps:
\begin{enumerate}
    \item dimension reduction;
    \item constructing the aPC basis;
    \item building the NN-aPC as a surrogate model of aPC modes and train the network for each mode.
\end{enumerate}
The trained NN-aPC can then be used to calculate the statistics of our QoI and to predict new instances of the continuous trajectories of the QoI, with newly collected sensor data. We will explain each of three steps and the prediction procedure below.

\subsubsection{Dimension reduction with principal component analysis}
As the first step, we find a lower dimensional random space spanned by a set of hidden random variables for the dimension reduction of our QoI. The most convenient way to do this is via the principal component analysis (PCA). Naturally, we would analyze the data of $k$, which is the source of randomness in Equation \ref{eqn:SDE}. Let $K$ be the $N_k\times N_k$ covariance matrix for the sensor measurements on $k$, i.e.,
\begin{equation}
    K_{i,j} = \Cov(k^{(i)}, k^{(j)}).
\end{equation}
Let $\lambda_l$ and $\phi_l$ be the $l$-th largest eigenvalue and its associated normalized eigenvector of $K$. Therefore, PCA yields
\begin{equation}\label{eqn:pca}
    K = \bm{\Phi}^T\bm{\Lambda}\bm{\Phi},
\end{equation}
where $\bm{\Phi} = [\phi_1, \phi_2, ..., \phi_{N_k}]$ is an orthonormal matrix and $\bm{\Lambda} = \text{diag}(\lambda_1, \lambda_2, ...\lambda_{N_k})$ is a diagonal matrix. Let $\bm{k}_s = [k_s^{(1)}, k_s^{(2)}, ..., k_s^{(N_k)}]^T$ be the results of the $k$ measurements of the $s$-th snapshot, then
\begin{equation}\label{eqn:get_xi}
    \bm{\xi}_s = \bm{\Phi}^T \sqrt{\bm{\Lambda}}^{-1} \bm{k}_s
\end{equation}
is an uncorrelated random vector, and hence $\bm{k}_s$ can be rewritten as a reduced dimensional expansion
\begin{equation}\label{eqn:dimension_reduction}
    \bm{k}_s \approx \bm{k}_0 + \sqrt{\bm{\Lambda}^M} \bm{\Phi}^M \bm{\xi}^M_s, \quad M < N_k,
\end{equation}
where $\bm{k}_0 = \mathbb{E}[\bm{k}]$ is the mean of each sensor's measurements. The choice of $M$ depends on how much energy should be maintained in the low-dimensional random space. For reasons of simplicity, we shall always use the reduced dimensional representation and omit the superscript $M$. We note that Equation \ref{eqn:dimension_reduction} can also be written in the form of the Karhunen-Lo\`eve expansion:
\begin{equation}\label{eqn:KL}
    k(x_k^{(i)};\omega_s) \approx k_0(x_k^{(i)}) + \sum_{l = 1}^M \sqrt{\lambda_l} k_l(x_k^{(i)}) \xi_{s,l}, \quad M < N_k,
\end{equation}
where $k_0(x_k^{(i)})$ is the mean of $k$ measurements at $x_k^{(i)}$, $k_l(x)$ is the $l$-th mode function of $k(x;\omega)$ whose value at $x_k^{(i)}$ coincides with the $i$-th entry of the eigenvector $\phi_l$, and $\xi_{s,l}$ is the $l$-th entry of the random vector $\bm{\xi}_s$. We want to extend the range of $k_l$ to the entire domain to approximate the continuous samples of $k(x;\omega)$. 

\subsubsection{Arbitrary polynomial chaos}
Assume that we have a set of $M$-dimensional samples of random vectors
\begin{equation*}
S:=\{\bm{\xi}_s\}_{s=1}^{N}
\end{equation*}
with hidden probability measure $\rho(\bm{\xi})$. Given a sufficiently large number of snapshots, we can approximate the underlying probability measure $\rho(\bm{\xi})$ by the discrete measure $\nu_S(\bm{\xi})$,
\begin{equation}\label{eqn:measure}
  \rho(\bm{\xi})\approx \nu_S(\bm{\xi})=\frac{1}{N}\sum_{\bm{\xi}_s\in S}\delta_{\bm{\xi}_s}(\bm{\xi}),
\end{equation}
where $\delta_{\bm{\xi}_s}$ is the Dirac measure. Then, a set of multivariate orthonormal polynomial basis $\{\psi_{\alpha}(\bm{\xi})\}_{\alpha=0}^P$ can be constructed via the Gram-Schmidt orthogonalization process following~\cite{Witteveen06, Leihuan2018}. The subscript $\alpha$ is the graded lexicographic multi-index and the number of basis, $P+1$, depends on the highest allowed polynomial order $r$ in $\psi_{\alpha}(\bm{\xi})$, following the formula
\begin{equation}
P + 1 = \frac{(r+M)!}{r!\,M!}.
\end{equation}
Specifically, the basis $\{\psi_{\alpha}(\bm{\xi})\}_{\alpha=0}^P$ are constructed using the recursive algorithm
\begin{equation}
  \psi_{\alpha}(\bm{\xi}) = w_{\alpha}^{\alpha} \psi^*_{\alpha}(\bm{\xi}) - \sum_{\beta \prec \alpha} w_{\beta}^{\alpha} \psi_{\beta} (\bm{\xi}),
\end{equation}
where $\psi^*_{\alpha}(\bm{\xi}):=\prod_{i=1}^M\bm{\xi}_i^{\alpha_i}$ represents the multivariate monomial basis function; the coefficients $w_{\beta}^{\alpha}$ are determined by imposing the orthonormal condition with respect to the discrete measure $\nu_S$, i.e.,
\begin{equation}
    \begin{aligned}
    \int\psi_{\alpha}(\bm{\xi}) \psi_{\beta}(\bm{\xi}) \dd \rho(\bm{\xi}) &\approx \int\psi_{\alpha}(\bm{\xi}) \psi_{\beta}(\bm{\xi}) \dd \nu_S(\bm{\xi})\\
    & = \frac{1}{N} \sum_{s=1}^{N} \psi_{\alpha}(\bm{\xi}_s) \psi_{\beta}(\bm{\xi}_s)\\
    & \equiv \delta_{\alpha,\beta}, \quad \beta \preceq \alpha.
    \end{aligned}
\end{equation}

With the polynomial basis $\{\psi_{\alpha}(\bm{\xi})\}$ that are automatically adapted to the distribution of $\bm{\xi}$, we can write any function $g(x;\bm{\xi})$ in the form of the aPC expansion,
\begin{equation}\label{eqn:aPC}
    g(x;\bm{\xi}) = \sum_{\alpha=0}^P g_{\alpha}(x) \psi_{\alpha}(\bm{\xi}),
\end{equation}
where the functions $g_{\alpha}(x)$ are called the aPC modes of $g$ and can be calculated by
\begin{equation}\label{eqn:calc_apc_modes}
    g_{\alpha}(x) = \frac{1}{N}\sum_{s=1}^{N}\psi_{\alpha}(\bm{\xi}_s)g(x;\bm{\xi}_s).
\end{equation}

\subsubsection{Learning stochastic modes}
\label{S:learning_modes}
\begin{figure}[htbp]
	\centering
	\includegraphics[width=0.6\linewidth]{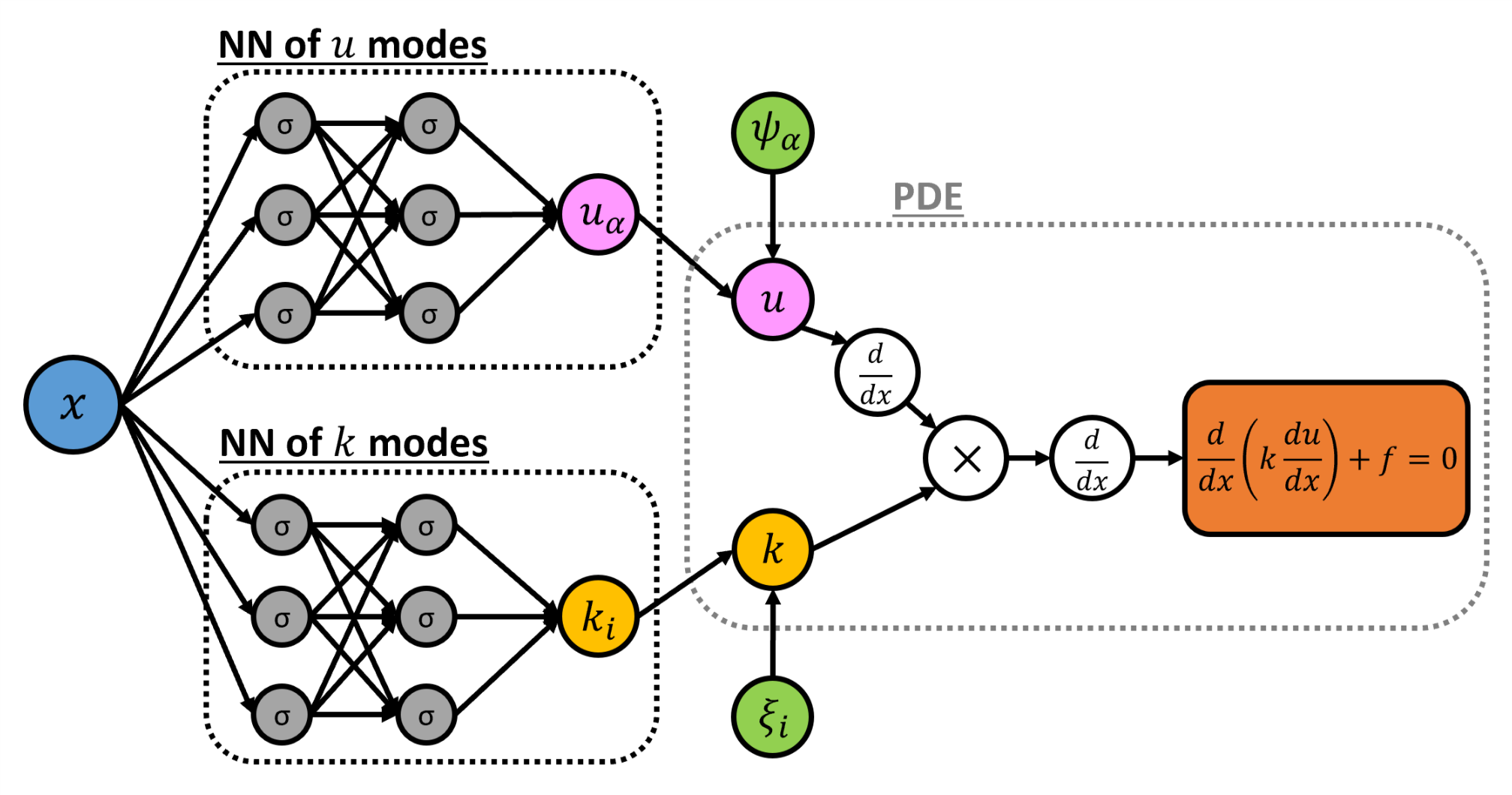}
	\caption{Schematic of the NN-aPC for solving the stochastic elliptic equation $-\dv{x}(k(x;\omega)\dv{x}u) = f$.}
	\label{fig:aPCiNN}
	\end{figure}

The key to our method is to train DNNs that predict the stochastic modes of our QoI. In this section, we focus on the inverse problem where we have to learn both the modes of $u$ and $k$. (Solving a forward problem is similar and more straightforward, and will be briefly discussed in Section \ref{S:4-1-1}.) Two disjoint DNNs are constructed, i.e., the network $\widehat{u_{\alpha}}$, which takes the coordinate $x$ as the input and outputs a $(P+1)\times 1$ vector of the aPC modes of $u$ evaluated at $x$, and the network $\widehat{k_i}$ that also takes the coordinate $x$ as the input and outputs a $(M+1)\times 1$ vector of the $k$ modes (we take $k_0$ in Equation \ref{eqn:KL} as the 0-th mode of $k$). Then, we can approximate $k$ and $u$ at the $s$-th snapshot by
\begin{equation}\label{eqn:recover_k}
    \tilde{k}(x; \omega_s) = \widehat{k_0}(x) + \sum_{i=1}^M \sqrt{\lambda_i} \widehat{k_i}(x) \xi_{s,i},
\end{equation}
and
\begin{equation}\label{eqn:recover_u}
    \tilde{u}(x; \omega_s) = \sum_{\alpha=0}^P \widehat{u_{\alpha}}(x) \psi_{\alpha} (\bm{\xi_s}).
\end{equation}
Similar to the PINN method, we construct the residual network via automatic differentiation and arithmetic operations of DNNs by substituting $u(x;\omega)$ and $k(x;\omega)$ in Equation \ref{eqn:SDE} with $\tilde{u}(x; \omega_s)$ and $\tilde{k}(x; \omega_s)$. This residual network is aPC-informed, because instead of being an uninterpretable black-box this network is designed to reflect the essence of the aPC expansion (see Figure \ref{fig:aPCiNN} for a schematic of the NN-aPC).

\begin{figure}[htbp]
	\centering
	\includegraphics[width=0.7\linewidth]{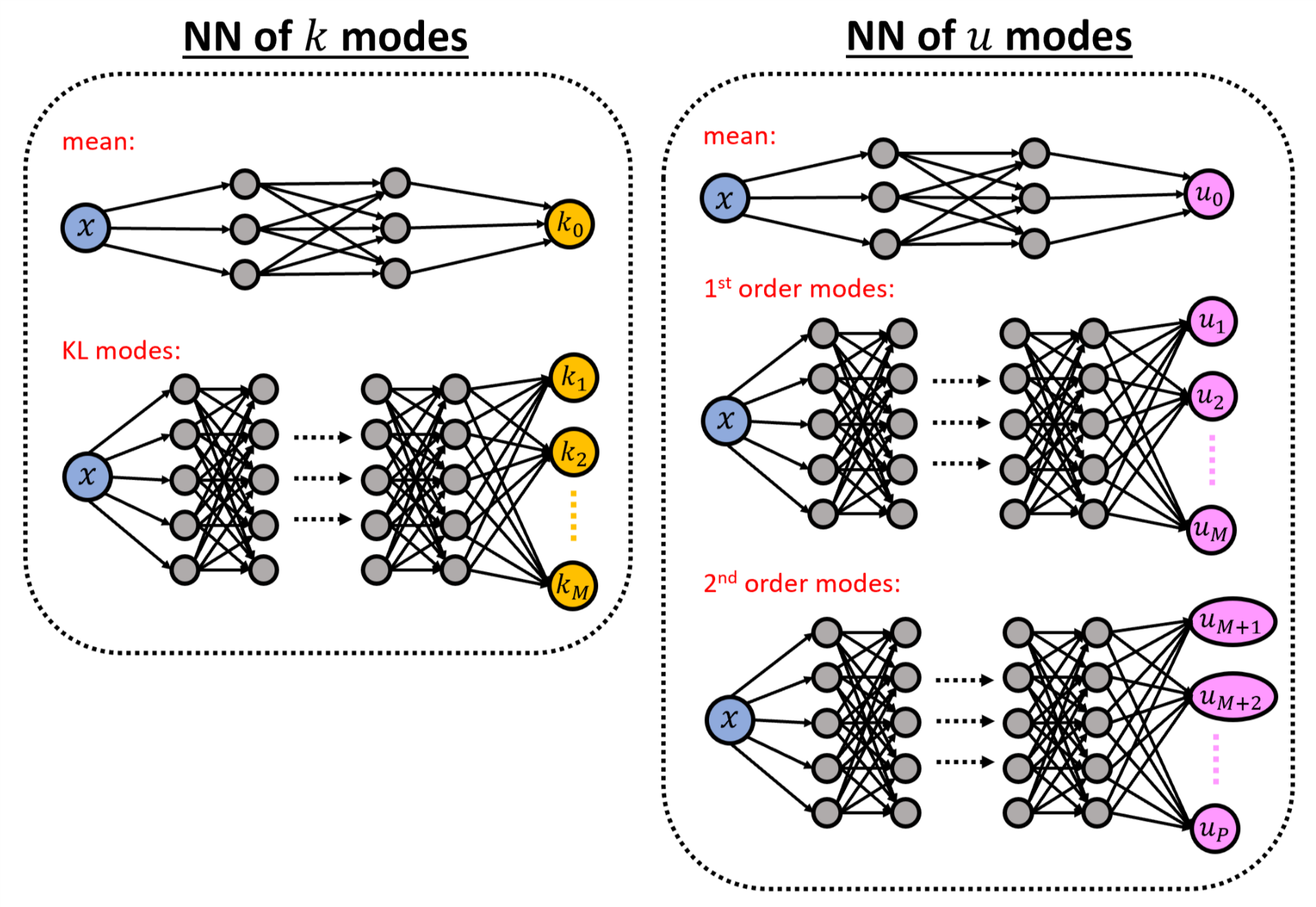}
	\caption{Schematic of the DNNs used for learning the stochastic modes of $k$ (left-hand side plot) and $u$ (right-hand side plot). The mean functions are modeled separately using small scale DNNs. For $k$, all its rest modal functions are modeled using one DNN. For $u$, the modes that correspond to the same order of aPC expansion are grouped together and are modeled with a single DNN.}
	\label{fig:nn_details}
	\end{figure}

In practice, for $\widehat{k_i}$, we separate the mean from the rest of the modes and learn it with a small scale DNN, and for $\widehat{u_{\alpha}}$, we group the modes corresponding to the same order of aPC expansion together and learn each group of modes with a separate DNN, as depicted in Figure \ref{fig:nn_details}. This is due to fact that the mean and the modes of different orders often correspond to vastly different scales.

The loss function is defined as a sum of the mean squared errors (MSE):
\begin{equation}\label{eqn:loss}
    \mathcal{L}(\mathcal{S}_t) = MSE_u + MSE_k + MSE_f,
\end{equation}
where
\begin{equation}
    MSE_u = \frac{1}{N N_u}\sum_{s=1}^N \sum_{i=1}^{N_u} \left[ \left( \tilde{u}(x_u^{(i)}; \omega_s) - u(x_u^{(i)};\omega_s) \right)^2 \right],
\end{equation}
\begin{equation}
    MSE_k = \frac{1}{N N_k}\sum_{s=1}^N \sum_{i=1}^{N_k} \left[ \left(\tilde{k}(x_k^{(i)}; \omega_s) - k(x_k^{(i)};\omega_s) \right)^2 \right],
\end{equation}
and
\begin{equation}
    MSE_f = \frac{1}{N N_f}\sum_{s=1}^N \sum_{i=1}^{N_f} \left[ \left(\mathcal{N}_x[ \tilde{u}(x_f^{(i)}; \omega_s) ; \tilde{k}(x_f^{(i)}; \omega_s) ] \right)^2 \right].
\end{equation}
So far, we have specified the training data, constructed the DNNs and formalized the loss function, and we now ready to train the DNNs.

\subsubsection{Predicting stochastic realizations}
\label{S:prediction}
In real applications, the training set could be generated from historical data, and the training process should be performed at the offline stage. The fine-tuned model shall be used to predict new random instances provided with new snapshots of sensor data, at the online stage. Suppose we have a snapshot of sensor data $\{\{(x_k^{(i)}, k_{\text{new}}^{(i)})\}_{i=1}^{N_k}, \{ (x_u^{(i)}, u_{\text{new}}^{(i)}) \}_{i=1}^{N_u} \}$; the first step is to extract the hidden random variables $\bm{\xi}_{\text{new}}$ from $\{ k_{\text{new}}^{(i)} \}_{i=1}^{N_k}$ using Equation \ref{eqn:get_xi}. Then, for any assigned location $x$, we can predict the modal functions for both $k(x;\omega)$ and $u(x;\omega)$ from the trained DNNs $\widehat{k_i}$ and $\widehat{u_{\alpha}}$. Finally, the prediction of $k$ and $u$ for the new random instance can be made via Equation \ref{eqn:recover_k} and \ref{eqn:recover_u}, respectively.

\subsection{Dropout for uncertainty}
\label{S:3-3}
Although DNNs can be used to approximate any measurable function accurately, standard DNNs do not capture model uncertainty. Dropout is one convenient way to quantify the approximation uncertainty in DNNs. The key idea of dropout is to drop units from the DNN independently and randomly with a pre-selected probability $p \in (0, 1)$. In the original work, dropout was only used during training, while no units were dropped at test time, i.e., the prediction of the DNN for the unknown data was deterministic. In the dropout for uncertainty, the units are also dropped at test time, resulting in stochastic predictions each time. The loss in the dropout inference is the summation of the original loss and a $l_2$ regularization term over the DNN parameters:
\begin{equation}\label{eqn:dropout_loss}
\mathcal{L} = \frac{1}{N_t}\sum_{i=1}^{N_t} l(\widehat{y_i}, y_i) + \lambda \sum_{i} \theta_i^2,
\end{equation}
where $N_t$ is the number of training points, $\widehat{y_i}$ is the prediction, $y_i$ is the true value, $l(\cdot, \cdot)$ is the loss for a single prediction, $\theta_i$ is any weight and bias, and $\lambda$ is the $l_2$ regularization rate. During prediction, the mean of the output is directly estimated by the Monte Carlo (MC) method,
\begin{equation}
\mathbb{E}(y) \approx \frac{1}{T}\sum_{t=1}^T \mathcal{NN}_t(x),
\end{equation}
where $\mathcal{NN}_t$ is the dropped neural network at the $t$-th prediction. The output variance is also estimated from these MC outputs.

\begin{figure}[htbp]
	\centering
	\includegraphics[width=0.5\linewidth]{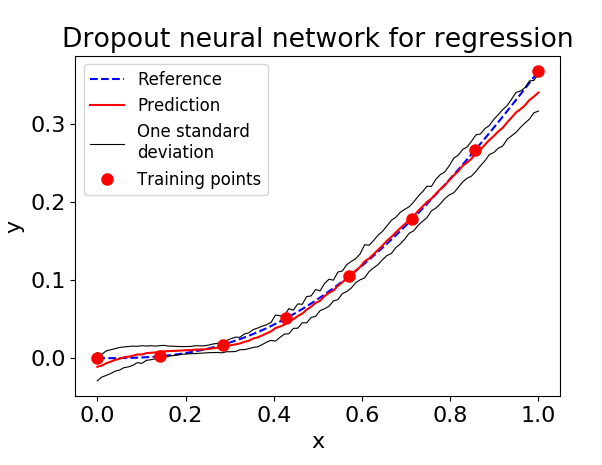}
	\caption{Dropout for uncertainty: An example of using the dropout in DNN to approximate the function $y = x^3e^{-x}$ in the domain $[0,1]$, where we use 4 hidden layers and 20 neurons per hidden layer, and we choose $p = 0.01, \lambda=10^{-6}$. The mean and standard deviation are calculated from 1000 MC samples.}
	\label{fig:dropout_regression}
	\end{figure}

Figure \ref{fig:dropout_regression} shows an example of using the dropout DNN for regression. We plan to use the dropout strategy in the NN-aPC method to estimate the uncertainty of our DNN model and as a guidance for active learning.

\section{Numerical Examples}
\label{S:4}
\subsection{Solving stochastic PDEs}
\label{S:4-1}
We first demonstrate the effectiveness of solving stochastic differential equations with the NN-aPC method for the forward and inverse problems.

\subsubsection{Forward problem: stochastic Poisson equation (a pedagogical example)}
\label{S:4-1-1}
Consider the following one-dimensional stochastic Poisson equation with homogeneous boundary conditions:
\begin{equation}\label{eqn:stochastic_poisson}
\begin{gathered}
    -\dv[2]{x} u = f(x;\omega), \quad x\in[-1, 1] \text{ and } \omega \in \Omega,\\
    u(-1) = u(1) = 0.
\end{gathered}
\end{equation}
Here $\Omega$ is the random space, the forcing term $f(x;\omega)\sim\mathcal{GP}(f_0(x), \Cov(x, x'))$ is a Gaussian random process with mean $f_0(x) = 10\sin(\pi x)$ and a squared exponential covariance function
\begin{equation}
\label{eqn:cov}
    \Cov(x, x') = \sigma^2\exp\left(-\frac{(x - x')^2}{l_c^2}\right),
\end{equation}
where the standard deviation $\sigma=1.0$ and the correlation length $l_c=0.5$.

In this and the following examples, all data are generated by the MC sampling method. Specifically, we sample $N=1000$ snapshots of continuous $f(x;\omega)$ trajectories $\{f_s=f(x;\omega_s)\}_{s=1}^{N}$ and extract from $\{f_s\}_{s=1}^N$ the values where the $N_f$ (virtual) sensors are located. For every $f(x;\omega)$ trajectory, we solve for its corresponding solution trajectories $u(x;\omega)$ using the finite difference method, and will use the statistics of these $u$ trajectories as our reference. To evaluate the performance of the trained model, we collect another $N_s=500$ snapshots of continuous $f(x;\omega)$ and $u(x;\omega)$ trajectories independently from the training data. The $N_s$ pairs of $(f, u)$ trajectories form our test sample set. Similarly, we extract the $f$-sensor data from every snapshot in the test set as the input at the predicting stage. We shall use the same $N$ and $N_s$ in the following tests, if not explicitly mentioned.

\begin{figure}[htbp]
	\centering
	\begin{subfigure}{0.45\textwidth}
		\centering
		\includegraphics[width=\linewidth]{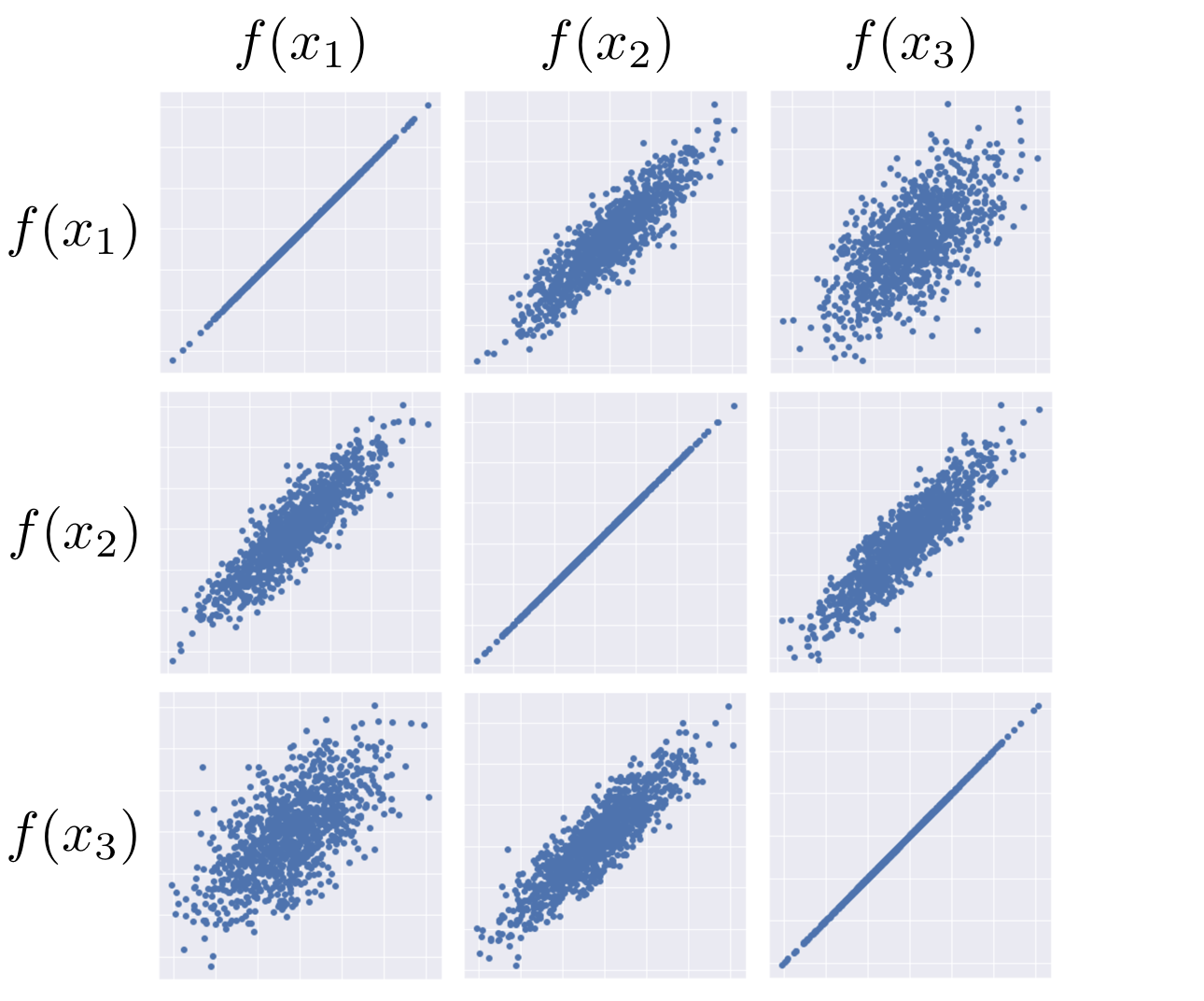}
		\caption{}
		\label{fig:correlation_sensor}
	\end{subfigure}%
	\begin{subfigure}{0.45\textwidth}
		\centering
		\includegraphics[width=\linewidth]{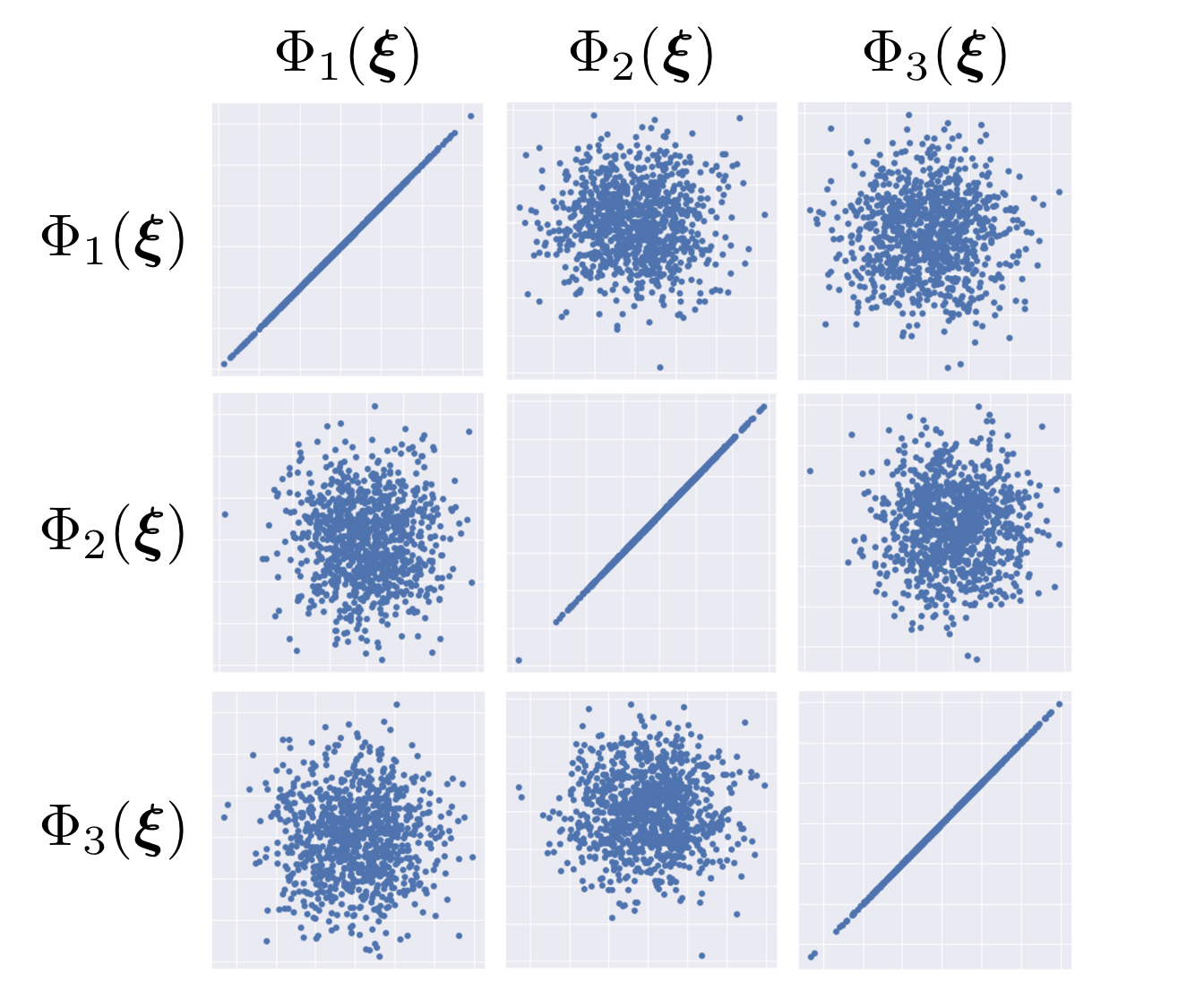}
		\caption{}
		\label{fig:correlation_mode}
	\end{subfigure}
	\caption{Correlation structure: (a) Scattered plots of data collected from the first three $f$-sensors. The correlation between different sensors is significant, and the closer the sensors are, the more correlated measurements they produces. (b) Scattered plots of the first three aPC basis evaluations. There is no correlation between different aPC basis functions.}
	\label{fig:correlation}
\end{figure}

We place $N_f=13$ sensors of $f(x;\omega)$ in the $[-1,1]$ domain (the sensors are equidistant) and keep $6$ principal random variables corresponding to $99\%$ stochastic energy after performing PCA. The solution $u(x;\omega)$ is approximated with a first-order aPC expansion. Figure \ref{fig:correlation} shows the scattered plots of the measurements from the first three $f$-sensors and the first three arbitrary polynomial basis. It is evident that the raw data from the measurements are correlated while their induced polynomials are not, thus the induced polynomials would serve as a valid set of basis in the random space.

\begin{figure}[htbp]
	\centering
	\begin{subfigure}{0.5\textwidth}
		\centering
		\includegraphics[width=\linewidth]{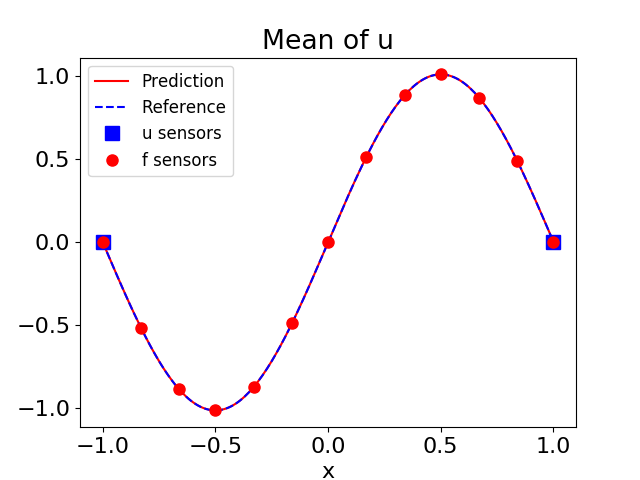}
		\caption{}
		\label{fig:411u_mean}
	\end{subfigure}%
	\begin{subfigure}{0.5\textwidth}
		\centering
		\includegraphics[width=\linewidth]{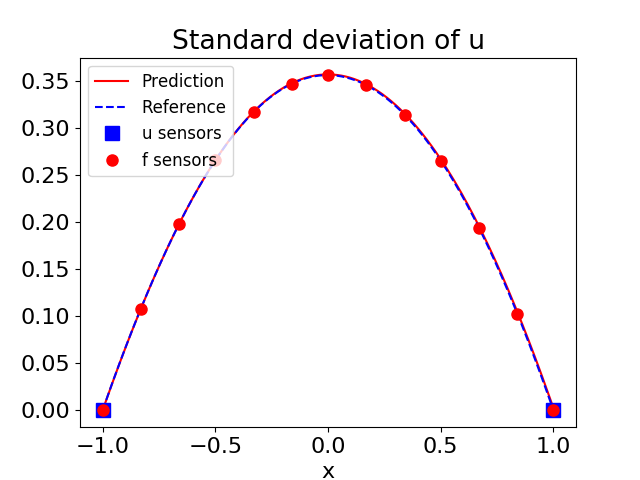}
		\caption{}
		\label{fig:411u_std}
	\end{subfigure}
	\caption{Forward problem: The mean function (a) and the standard deviation (b) of $u$ predicted using the trained DNN. The reference is calculated from the 1000 snapshots of continuous $u$ samples. In this case, 13 $f$-sensors are employed, while only 2 $u$-sensors are placed at the domain boundaries.}
	\label{fig:411u}
\end{figure}

\begin{figure}[htbp]
	\centering
	\includegraphics[width=\linewidth]{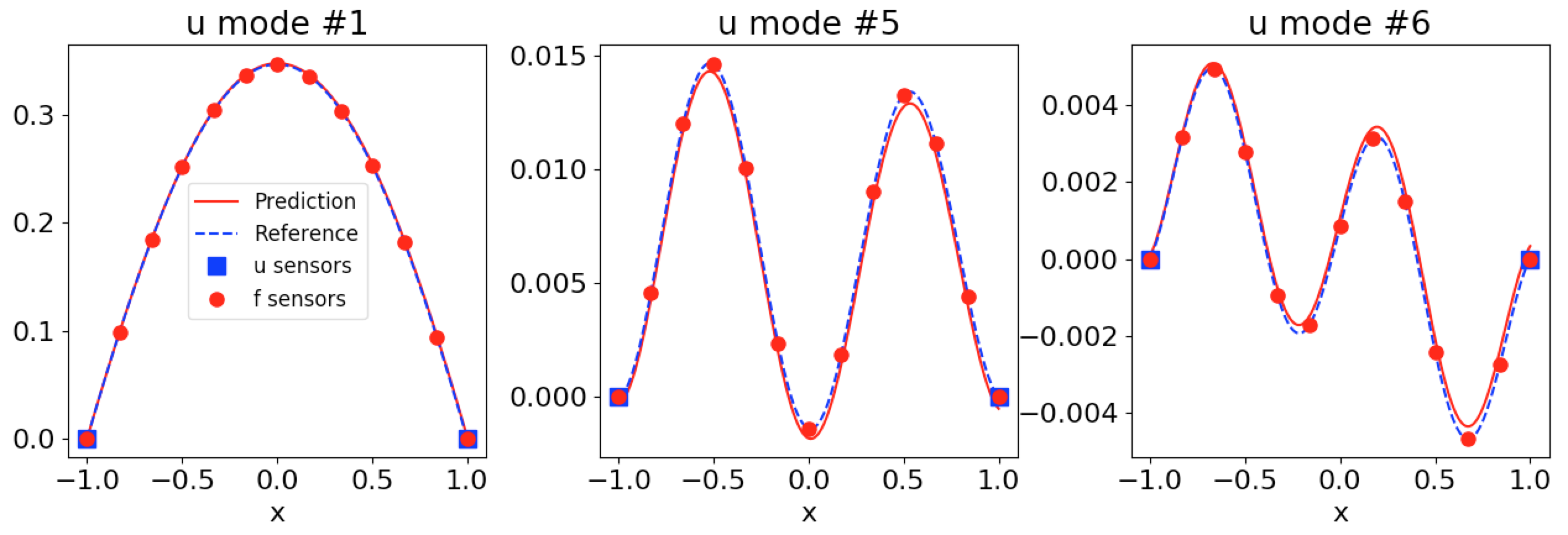}
	\caption{Forward problem: Plots of three (out of six) aPC modes of $u$ versus the reference, which is calculated from Equation \ref{eqn:calc_apc_modes} using the 1000 continuous $u$ samples. The predicted modes match the true modes closely.}
	\label{fig:411u_modes}
\end{figure}

The DNNs used to approximate the modes of $u$ are constructed as in Figure \ref{fig:nn_details}, where we use an isolated small scale DNN of 2 hidden layers with 4 neurons per hidden layer to approximate the mean profile, and a DNN of 4 hidden layers with 32 neurons per hidden layer to model the modes. The {\it tanh} function is selected as the default activation function due to it is second order differentiable. Then, a DNN for the residual can be constructed via auto-differentiation and arithmetic operations. The training set $\mathcal{S}_t$ is
\begin{equation*}
    \mathcal{S}_t = \left\{\{(-1, 0), (1, 0)\}, \{(x_f^{(i)}, 0\}_{i=1}^{N_f} \right\}_{s=1}^N,
\end{equation*}
where the $u$ data is collected only at the boundaries to provide boundary conditions. The loss function is slightly modified based on Equation \ref{eqn:loss} to add a $l_2$ regularization term. At the training stage, we choose the $l_2$ regularization rate $\lambda=0.001$, and use the Adam~\cite{Kingma15} optimizer with learning rate $0.001$ to train our model for $20000$ epochs. Figure \ref{fig:411u} shows the predicted mean and standard deviation of the solution $u$ versus the reference. Figure \ref{fig:411u_modes} shows our DNN prediction of three $u$ modes where the reference modes are calculated by Equation \ref{eqn:calc_apc_modes}. We can see that the NN-aPC method makes accurate predictions of the mean and standard deviation of the solution $u(x;\omega)$ and learns the arbitrary polynomial chaos modes. 

\begin{figure}[htbp]
	\centering
	\begin{subfigure}{0.5\textwidth}
		\centering
		\includegraphics[width=\linewidth]{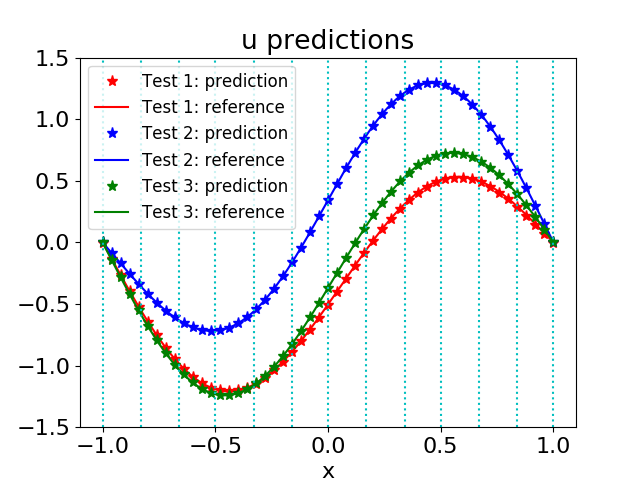}
		\caption{}
		\label{fig:411_prediction}
	\end{subfigure}%
	\begin{subfigure}{0.5\textwidth}
		\centering
		\includegraphics[width=\linewidth]{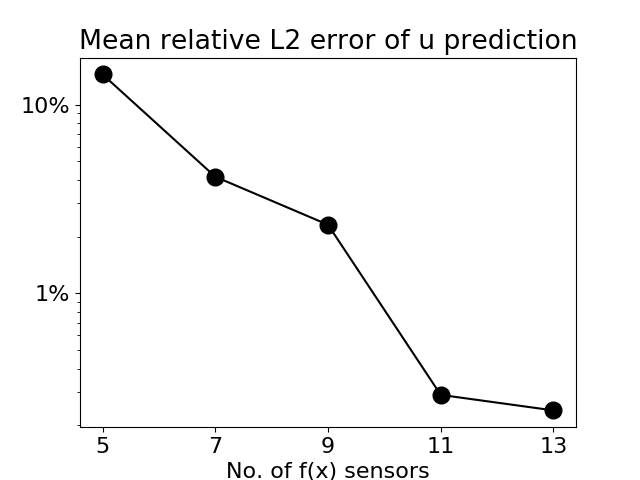}
		\caption{}
		\label{fig:411_error}
	\end{subfigure}
	\caption{Making predictions for the forward problem using the NN-aPC: (a) For three different snapshots in the test data set, we compare the predicted solution $u$, calculated using the measurements of 13 $f$-sensors whose locations are marked with the dashed lines, versus the true $u$. (b) Relative $L_2$ error of the predicted $u$, averaged for all snapshots in the test data set, versus the number of $f$-sensors placed in the domain.}
\end{figure}

The trained model is then used to predict the QoI at any location $x_o$ given new snapshots of sensor data, with only one forward evaluation of the DNN with $x_o$ as the input, as described in Section \ref{S:prediction}. Figure \ref{fig:411_prediction} illustrates the prediction of solution for three different snapshots of the sensor data in the test samples; our prediction recovers the true solution very well. In Figure \ref{fig:411_error}, we use an increasing value of $N_f$ to study the effect of the number of $f$-sensors on the accuracy of the predicted solutions; we observe that when more $f$-sensors are deployed, we can achieve better accuracy.

\subsubsection{Inverse problem: stochastic elliptic equation}
\label{S:4-1-2}
We solve the one-dimensional stochastic elliptic equation as an inverse problem, where we have some extra information on the solution $u(x;\omega)$ but incomplete information of the diffusion coefficient $k(x;\omega)$. The equation reads
\begin{equation}
\label{eqn:stochastic_elliptic}
\begin{gathered}
    -\dv{x}\left(k(x;\omega)\dv{x}u\right) = f(x), \quad x\in[-1, 1] \text{ and } \omega \in \Omega,\\
    u(-1) = u(1) = 0.
\end{gathered}
\end{equation}
In this example, we use a constant forcing term $f(x)=10$. The randomness comes from the diffusion coefficient $k(x;\omega)$, for which we only have limited information at the locations where we place the $k$-sensors. Here in this example, $k$ is modeled and sampled from a non-Gaussian random process such that
\begin{equation}
    \log(k(x;\omega)) \sim\mathcal{GP}(k_0(x),\Cov(x,x')),
\end{equation}
where the mean $k_0(x) = \sin(3\pi x/2)/5$ and the covariance function has the same form as in Equation \ref{eqn:cov}, where we set the standard deviation $\sigma=0.1$ and correlation length $l_c=1.0$. We use the same strategy as in Section \ref{S:4-1-1} to generate the training and testing samples, and the training set is constructed as in Equation \ref{eqn:training_set}.

\begin{figure}[htbp]
	\centering
	\begin{subfigure}{0.5\textwidth}
		\centering
		\includegraphics[width=\linewidth]{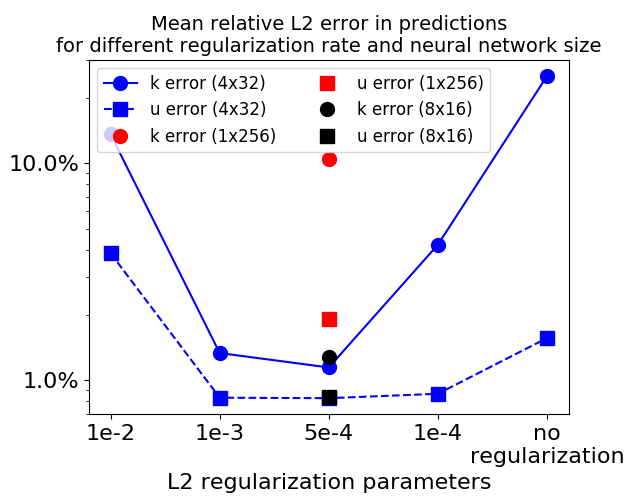}
		\caption{}
		\label{fig:412_errorvslambda}
	\end{subfigure}%
	\begin{subfigure}{0.5\textwidth}
		\centering
		\includegraphics[width=\linewidth]{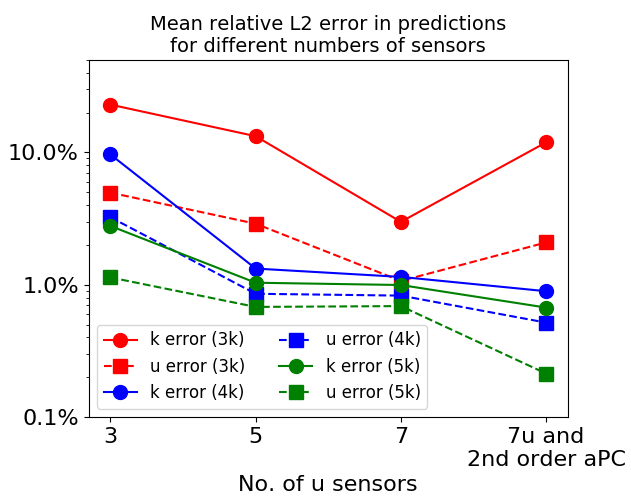}
		\caption{}
		\label{fig:412_errorvsnsensor}
	\end{subfigure}
	\caption{Comparing prediction accuracy for the inverse problem using the NN-aPC: (a) The mean of relative $L_2$ error in predicting $u$ and $k$ trajectories in the test set when using different sizes of DNNs and different $l_2$ regularization rate. In this case, we use 1st-order aPC expansion 4 $k$-sensors and 7 $u$-sensors are deployed and both DNNs have the same size. (b) The mean relative $L_2$ error in predicting $u$- and $k$-trajectories versus the number of sensors deployed. In this case we use the 1st-order aPC expansion, choose $\lambda = 0.0005$, and the DNNs have 4 hidden layers with 32 neurons per hidden layer.}
\end{figure}

For this case, two groups of DNNs, i.e., $\widehat{k_i}$ and $\widehat{u_{\alpha}}$ are built to calculate the modes for $k$ and $u$. Again, we use the Adam optimizer with learning rate 0.001 to train the DNNs for 50000 epochs. In Figure \ref{fig:412_errorvslambda}, we use the 1st-order aPC expansion and we study the impact of using different $l_2$ regularization rate $\lambda$ and different shapes of DNNs. We only change the DNNs that learn the stochastic modes, while the DNNs that learn the mean profiles are fixed to have 2 hidden layers and 4 neurons per hidden layer; this is the default setting for the future examples as well. In these plots, we compare the averaged relative $L_2$ error of predicting $k$ and $u$ in the test set. The results indicate that a moderate choice of $\lambda$ (0.0005) gives us the most accurate predictions, since a too small/large choice of $\lambda$ causes over-/under- fitting. A suitable choice of DNN shape (4 hidden layers with 32 neurons per hidden layer) produces the best trained model. For the rest of this numerical example, we shall adopt these optimal DNN setting.

In Figure \ref{fig:412_errorvsnsensor}, we use the 1st-order aPC expansion and compare the averaged relative $L_2$ error in predictions when different numbers of $k$- and $u$-sensors are deployed to collect the training data. In general, the proposed method makes more accurate predictions after training with data collected from more $k$- and $u$-sensors. One reason is that a larger number of sensors supports a greater variety of input $x$ in the training data, thus feeding more information to the model to reduce the probability of over-fitting. Also, more $k$-sensors allows for a higher effective random dimension of the aPC expansion for better approximation. Figure \ref{fig:412_errorvsnsensor} also shows that a 2nd-order aPC expansion helps to improve predictions. We note that this is not the case when we use only three $k$-sensors. The bottleneck here is insufficient random dimension, so without enough training information, adopting the 2nd-order aPC expansion doubles the number of $\widehat{u_{\alpha}}$ net outputs and would only increase the risk of over-fitting.

\begin{figure}[htbp]
	\centering
	\begin{subfigure}{0.8\textwidth}
		\centering
		\includegraphics[width=\linewidth]{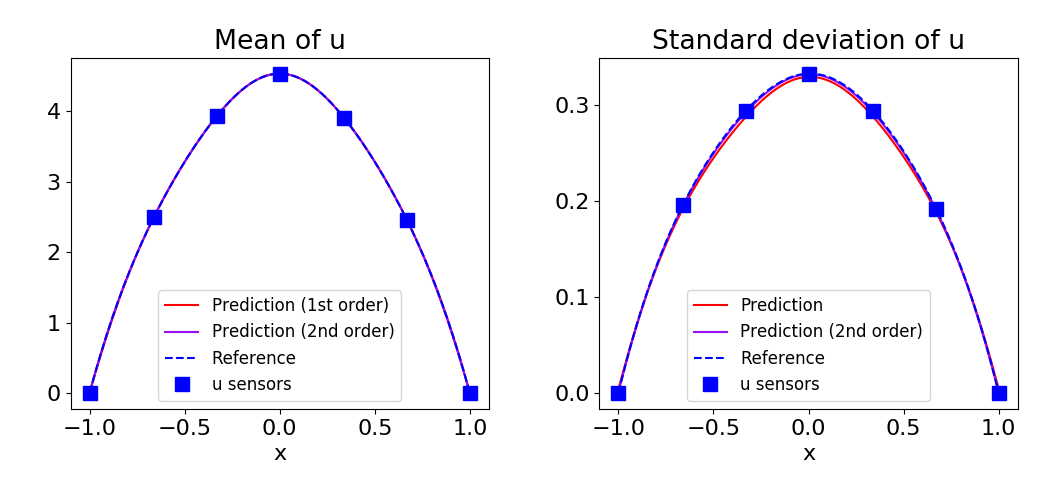}
		\caption{}
		\label{fig:412_umean_std}
	\end{subfigure}
	\begin{subfigure}{\textwidth}
		\centering
		\includegraphics[width=\linewidth]{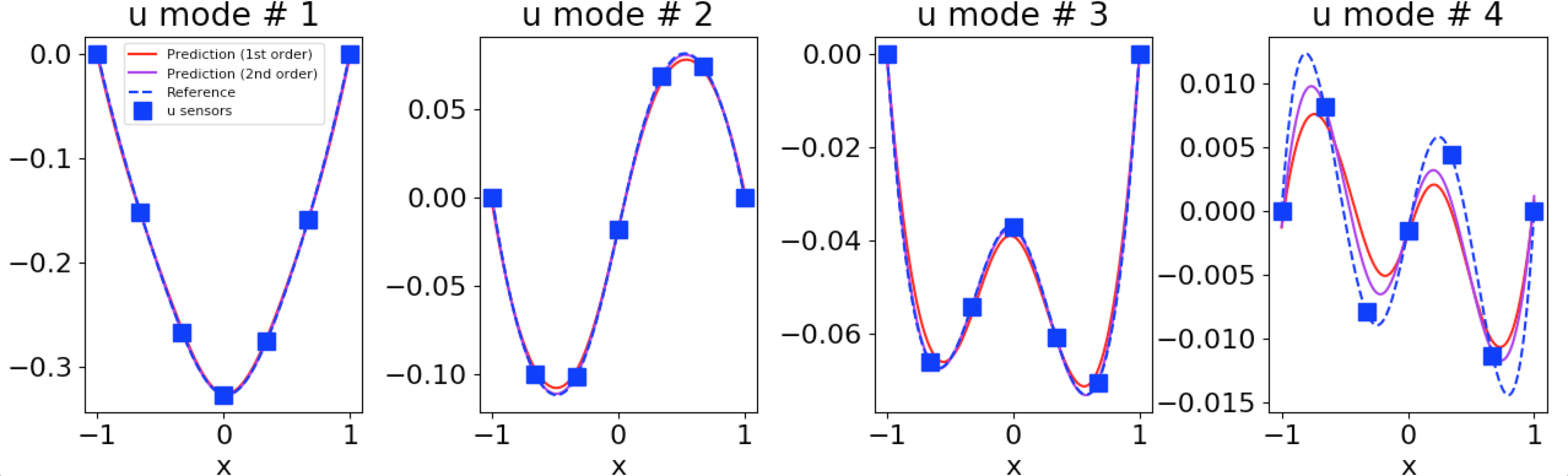}
		\caption{}
		\label{fig:412_umodes}
	\end{subfigure}
	\caption{Testing 2nd-order aPC expansions: (a) The predicted mean/standard deviation of $u$ calculated with a 1st- and 2nd-order aPC expansions versus the reference. (b) The first 4 modes of $u$ calculated by the 1st- and 2nd-order aPC expansion. The reference solutions in all plots are calculated with the continuous trajectories that produced the training data. The results are generated with 7 $u$-sensors (blue squares) and 4 $k$-sensors (red dots in Figure \ref{fig:412_k}).}
	\label{fig:412_u}
\end{figure}

\begin{figure}[htbp]
	\centering
	\begin{subfigure}{0.8\textwidth}
		\centering
		\includegraphics[width=\linewidth]{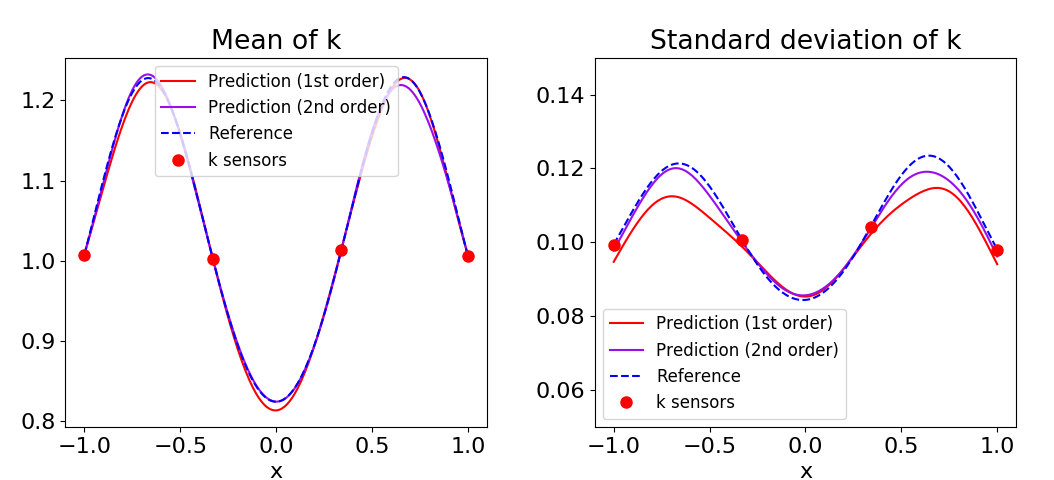}
		\caption{}
		\label{fig:412_kmean_std}
	\end{subfigure}
	\begin{subfigure}{\textwidth}
		\centering
		\includegraphics[width=\linewidth]{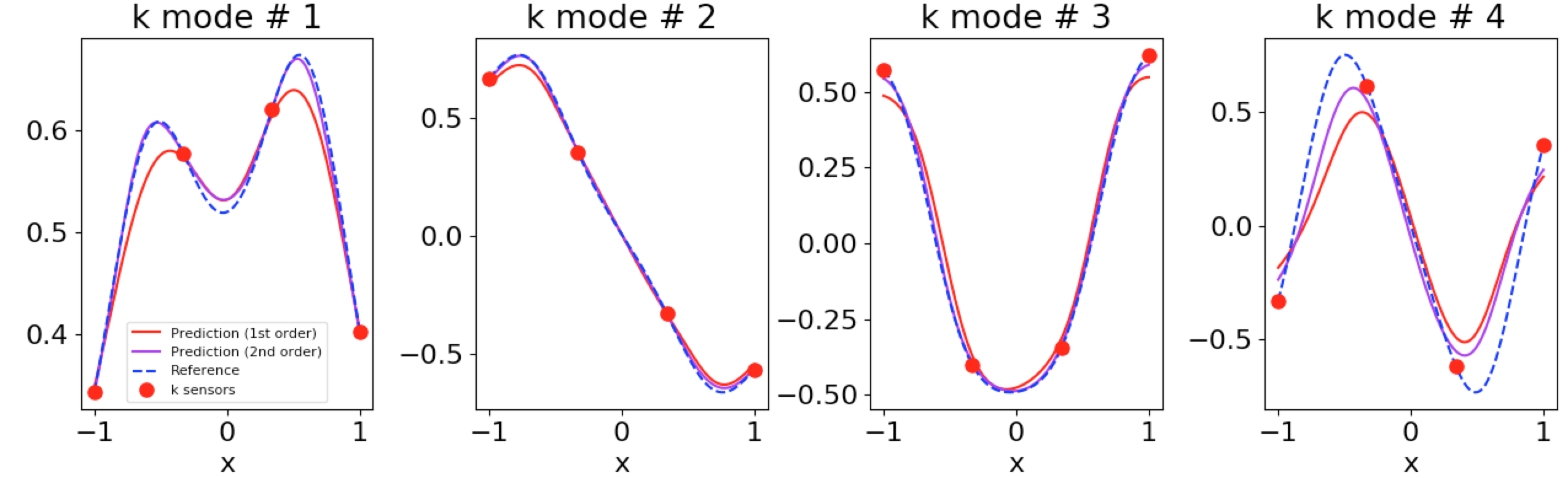}
		\caption{}
		\label{fig:412_kmodes}
	\end{subfigure}
	\caption{Testing 2nd-order aPC expansions: (a) The predicted mean/standard deviation of $k$ calculated with a 1st- and 2nd-order aPC expansions versus the reference. (b) The first 4 modes of $k$ calculated by the 1st- and 2nd-order aPC expansion. The reference in all plots are calculated with the continuous trajectories that produced the training data. The results are generated with the same setup of sensors as that of Figure \ref{fig:412_u}.}
	\label{fig:412_k}
\end{figure}

\begin{table}[htbp]
\centering
\begin{tabular}{|c|c|c|c|c|c|c|c|}
\hline
\multicolumn{1}{|l|}{} &  & mean & std & mode 1 & mode 2 & mode 3 & mode 4 \\ \hline
\multirow{2}{*}{$k$} & 1st-order & 0.54\% & 5.26\% & 4.15\% & 5.39\% & 12.03\% & 42.81\% \\ \cline{2-8} 
                     & 2nd-order & 0.45\% & 1.87\% & 1.28\% & 1.95\% &  3.67\% & 29.95\% \\ \hhline{|=|=|=|=|=|=|=|=|}
\multirow{2}{*}{$u$} & 1st-order & 0.14\% & 1.83\% & 0.98\% & 3.60\% &  4.34\% & 45.56\% \\ \cline{2-8} 
                     & 2nd-order & 0.14\% & 0.51\% & 0.08\% & 0.80\% &  1.04\% & 32.16\% \\ \hline
\end{tabular}
\caption{Comparing the relative $L_2$ error when using the 1st- and 2nd-order aPC expansion, and using data from 4 $k$-sensors and 7 $u$-sensors for training.}
\label{tab:412_error}
\end{table}

\begin{figure}[htbp]
	\centering
	\includegraphics[width=0.8\linewidth]{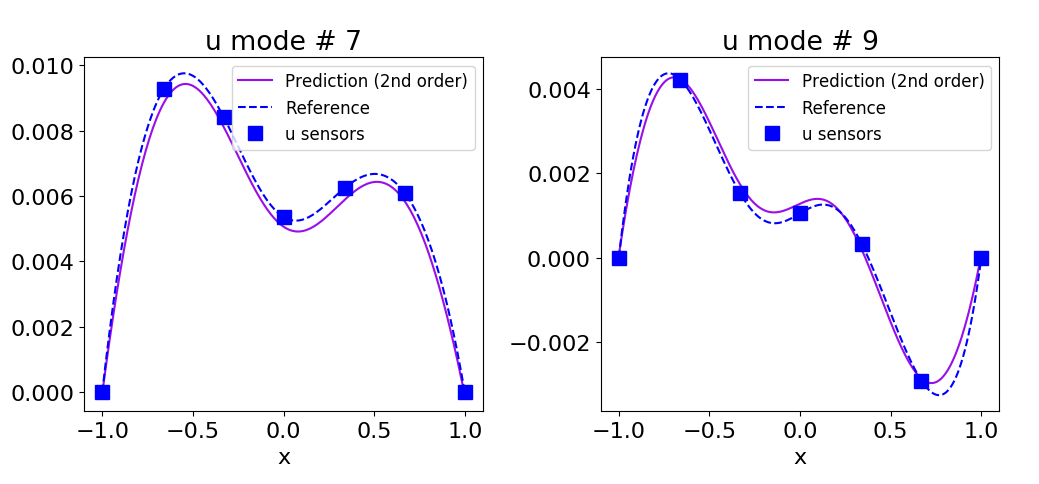}
	\caption{Testing 2nd-order aPC expansions: Two $u$-modes of higher order calculated with the 2nd-order aPC expansion. The magnitude of the higher-order modes is smaller compared to the lower-order modes, but the NN-aPC method is still able to capture the small magnitude modes. The results are generated with the same setup of sensors as that of Figure \ref{fig:412_u}.}
	\label{fig:412_highmode}
\end{figure}

We use a combination of 7 $u$-sensors and 4 $k$-sensors. Figure \ref{fig:412_umean_std} and \ref{fig:412_kmean_std} compare the mean and standard deviation of $u$ and $k$ calculated by the trained DNNs when we use the 1st- and 2nd-order aPC expansions. Figure \ref{fig:412_umodes} shows the first four common aPC modes of $u$ for both expansions, where we can see that the 2nd-order aPC expansion helps to improve the accuracy of the predicted lower-order modes. The 2nd-order aPC expansion would also help with learning the $k$ modes, as depicted in Figure \ref{fig:412_kmodes}, even if the $\widehat{k_i}$ network is disjoint with the $\widehat{u_{\alpha}}$ network. Table \ref{tab:412_error} lists the relative $L_2$ error of the trained model in calculating the mean, standard deviation and all the common modes of $k$ and $u$. We can see that using a 2nd-order aPC significantly improves the accuracy. We also note that in Figure \ref{fig:412_kmean_std}, due to the placement of $k$-sensors, the measurements from each sensor will yield almost the same mean (1.0) and standard deviation (0.1), but the trained $k$ net would reveal the non-trivial wavy structure of its mean and standard deviation in the entire domain, containing information more than just the training data could provide. This indicates that there is information fusion of all three types of training data and the stochastic differential equation, rather than a simple interpolation. Figure \ref{fig:412_highmode} shows that the higher-order aPC modes can also be captured even if they exist in a much smaller scale compared to the lower-order more energetic modes.

\begin{figure}[htbp]
	\centering
	\begin{subfigure}{0.5\textwidth}
		\centering
		\includegraphics[width=\linewidth]{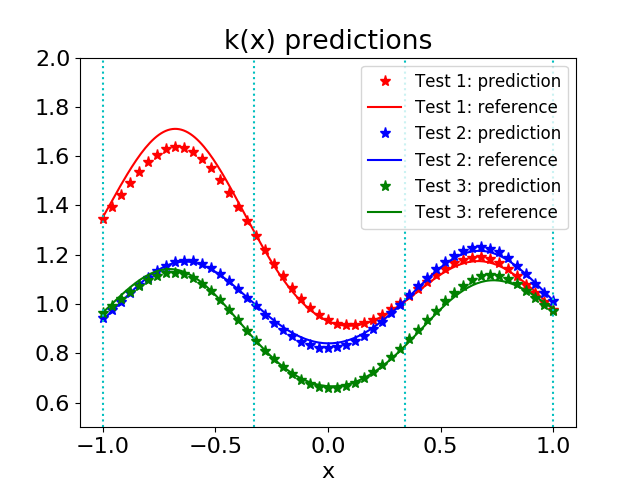}
		\caption{}
		\label{fig:412_kpred}
	\end{subfigure}%
	\begin{subfigure}{0.5\textwidth}
		\centering
		\includegraphics[width=\linewidth]{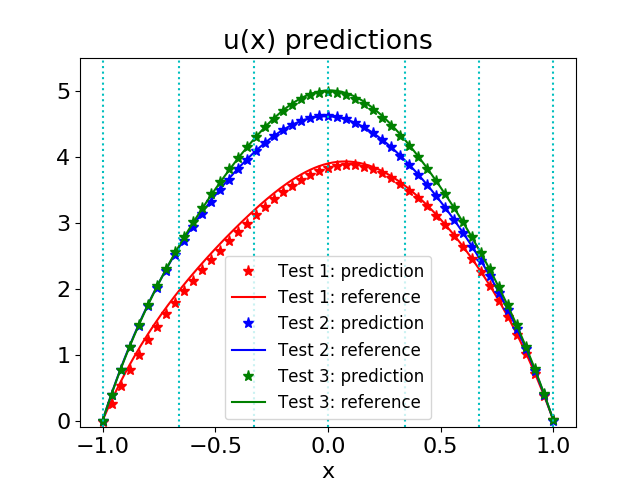}
		\caption{}
		\label{fig:412_upred}
	\end{subfigure}
	\caption{Making predictions for the inverse problem using the NN-aPC: For three different snapshots in the test data set, we compare the predicted solution $k$ (in plot (a)) and $u$ (in plot (b)), calculated using the measurements of 7 $u$-sensors and 4 $k$-sensors (denoted by the dashed lines).}
	\label{fig:412_pred}
\end{figure}

Figure \ref{fig:412_pred} shows the prediction of $k$ and $u$ for three arbitrary snapshots in the test sample set. Every pair of predictions are made based on a single time measurement from 4 $k$- and 7 $u$-sensors, and they agree with the true reference value. Again, the predicting stage takes little computational time since only one forward evaluation of the well-trained DNNs is needed for every input $x$, no matter how many snapshots of predictions are to be made. 

\subsection{Solving deterministic differential equations with PINN and dropout}
\label{S:4-2}
In this part, we first show that dropout reduces over-fitting in solving differential equations. Then, more importantly, we show that the dropout-induced uncertainty serves as useful guidance for active learning.

\subsubsection{Reducing over-fitting}
To show how dropout reduces over-fitting, we implement the dropout neural networks to solve both a forward Poisson equation and an inverse elliptic equation. They are the deterministic version of Equation \ref{eqn:stochastic_poisson} and \ref{eqn:stochastic_elliptic}. For the forward Poisson equation, we choose $f(x) = 9\pi^2\sin(3\pi x/2)/4$ as the forcing term. A dropout neural network with 6 hidden layers and 100 neurons per hidden layer is constructed to model the solution $u(x)$. The training data consists of 2 $u$ measurements at both boundaries and 6 $f$-sensors in the domain. For the inverse elliptic equation, we choose, as before, $f(x) = 10$ and $k(x) = \exp(\sin(3\pi x/2)/5)$ as the hidden diffusion coefficient, and we use 5 $k$-sensors and 7 $u$-sensors. Two separate DNNs are constructed: a small scale regular DNN that has 2 hidden layers with 4 neurons per hidden layer to model the function $u(x)$, and a large scale dropout neural network with 6 hidden layers and 100 neurons per hidden layer to model the function $k(x)$. The dropout rate in both examples is fixed at 0.01.

\begin{figure}[htbp]
	\centering
	\begin{subfigure}{0.5\textwidth}
		\centering
		\includegraphics[width=\linewidth]{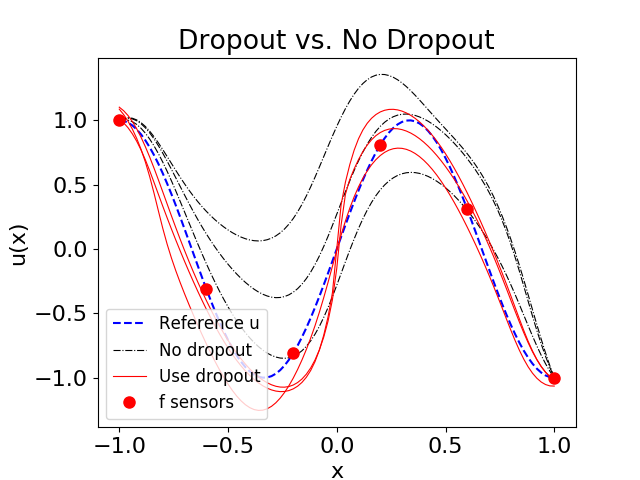}
		\caption{}
		\label{fig:421_dropoutforoverfitting}
	\end{subfigure}%
	\begin{subfigure}{0.5\textwidth}
		\centering
		\includegraphics[width=\linewidth]{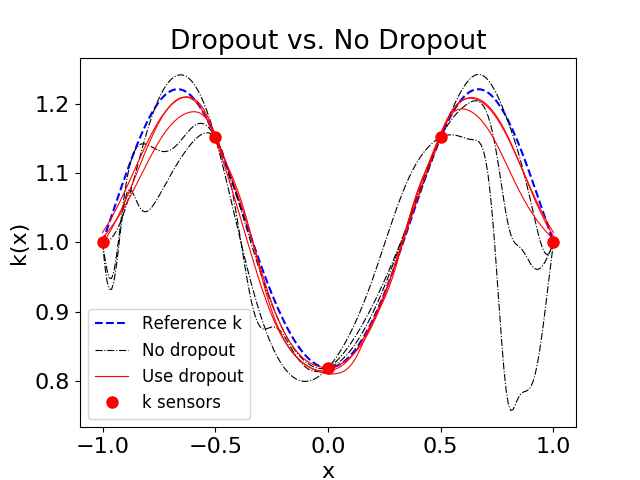}
		\caption{}
		\label{fig:422_dropoutforoverfitting}
	\end{subfigure}
	\caption{Dropout to reduce over-fitting: A comparison of the predicted QoI using the PINN/dropout and the regular PINN. For each case, three independent runs are conducted. (a) Forward problem: we solve the Poisson equation with 6 $f$-sensors (red dots) and 2 $u$-sensors (at the domain boundary, not shown in the plot). (b) Inverse problem: we solve an elliptic equation with 5 $k$-sensors (red dots) and 7 $u$-sensors (equidistantly placed in the domain, not shown in the plot).}
	\label{fig:42_overfitting}
\end{figure}

For both the forward and inverse problem, we train the DNNs as described in Section \ref{S:3-3}, using an Adam optimizer with learning rate $0.001$ for 30000 epochs. Due to the lack of sufficient number of sensors, there is a big chance that over-fitting occurs at the training stage if we use just regular DNNs. Figure \ref{fig:42_overfitting} shows a comparison of the results when we train the networks with and without using dropout. As we can see in the plots, the results from a regular DNN are very different from each other, showing irregular jumps of large amplitudes, while the results from dropout DNNs are similar to each other, and they are closer to the truth. This shows that dropout works as an effective means of reducing over-fitting.

\subsubsection{Estimating DNN approximation uncertainty}
More importantly, the uncertainty introduced by dropout serves as a useful guidance for active learning. Again, we solve the same inverse elliptic equation but this time we are provided with additional $k$-sensors. At the training stage, we add an $l_2$ regularization term with $\lambda=10^{-6}$ to the loss function. At the predicting stage, we evaluate the dropout network for $10000$ times to estimate the mean and standard deviation. Next, a new $k$-sensor is placed where the standard deviation reaches its maximum. If it happens that the location to add the new sensor is close to an existing $k$-sensor by a threshold distance of $\rho$, we do not add the new sensor, but instead, we count the nearest existing sensor twice as if we have added a virtual sensor at the same location of the existing one. In practice, here we choose $\rho = 0.03$. We have designed an automatic iterative procedure for active learning that reduces the cost of adding new sensors by efficiently re-using the old ones.

\begin{figure}[htbp]
	\centering
	\begin{subfigure}{0.5\textwidth}
		\centering
		\includegraphics[width=\linewidth]{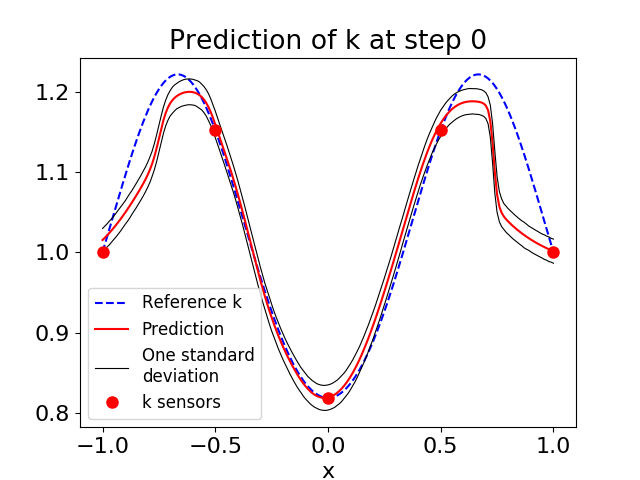}
		\caption{}
		\label{fig:422_first}
	\end{subfigure}%
	\begin{subfigure}{0.5\textwidth}
		\centering
		\includegraphics[width=\linewidth]{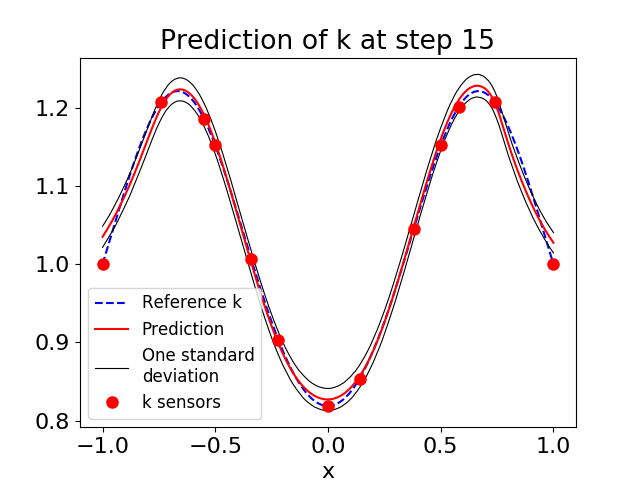}
		\caption{}
		\label{fig:422_last}
	\end{subfigure}
	\caption{Active learning using dropout approximation uncertainty: (a) The prediction of $k$ and its dropout-induced uncertainty of the starting step with 5 $k$-sensors. (b) The prediction of $k$ and its dropout-induced uncertainty of the last step with 13 $k$-sensors.}
	\label{fig:422_first_and_last}
\end{figure}

\begin{figure}[htbp]
	\centering
	\includegraphics[width=0.5\linewidth]{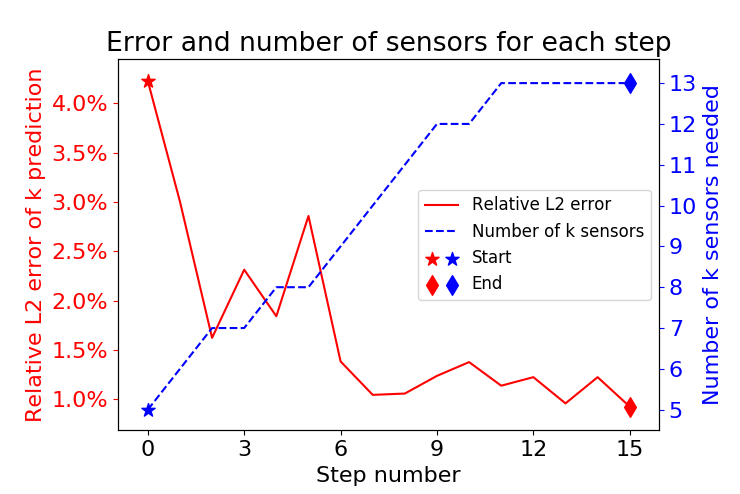}
	\caption{Effectiveness of active learning: The red solid line shows the relative $L_2$ error of the predicted $k$ revealing a decaying trend. The blue dashed line shows the number of $k$-sensors deployed in each step.}
	\label{fig:422_error}
\end{figure}

Figure \ref{fig:422_first_and_last} shows the initial prediction of $k$ and the prediction after iterating the above algorithm for 15 steps. In Figure \ref{fig:422_last}, the sensors are clustered where the curvature of $k$ is big, which is consistent with our intuition that we should put more sensors where the function changes rapidly. In Figure \ref{fig:422_error}, we see that during these 15 iterations, only 8 new sensors are deployed while the relative error of $k$ predictions reduces from more than $5\%$ to less than $1\%$.

\subsection{Active learning for inverse stochastic elliptic problems}
We consider the inverse stochastic elliptic problem as described in Equation \ref{eqn:stochastic_elliptic} for active learning, but in this example, $\log(k(x;\omega))$ is modeled by a Gaussian random process with correlation length $l_c=0.5$. To start with, we have 1000 snapshots of data from 3 $k$-sensors, 7 $u$-sensors and 21 $f$-sensors that are equidistantly distributed in the physical domain, and our goal is to infer $k(x;\omega)$ in the entire domain. Suppose we are then provided with additional sensors of $k$; we shall allocate them according to the uncertainty induced by the dropout neural networks. In practice, we model the modes of $k(x;\omega)$ with a dropout neural network that has 4 hidden layers with 128 neurons per hidden layer, and a dropout rate 0.01. The solution $u(x;\omega)$ is expanded with the 1st-order aPC expansion, while the modes are modeled by a regular DNN with 4 layers and 32 neurons per hidden layer. The reasons for implementing dropout only on $\widehat{k_i}$ are as follows:
\begin{itemize}
    \item[1.] Dropout can be used to efficiently reduce over-fitting, and the over-fitting issue of $k(x,\omega)$ is worse than that of $u(x,\omega)$.
    \item[2.] As an inverse problem, our main goal is to identify $k(x,\omega)$ and then identify where we should add more $k$-sensors to enhance the accuracy of prediction.
\end{itemize}
We use an Adam optimizer with learning rate 0.0005 to train the networks for 50000 epochs. The mean and standard deviation of the $k$ modes model are evaluated from 10000 independent evaluations of the dropout neural network. We place the additional $k$-sensor where the standard deviation of the first mode attains its maximum. This is because the first mode is associated with the largest eigenvalue, therefore bringing the largest impact to the stochastic structure, thus it is most important to learn the first mode accurately.

\begin{figure}[htbp]
	\centering
	\begin{subfigure}{0.33\textwidth}
		\centering
		\includegraphics[width=\linewidth]{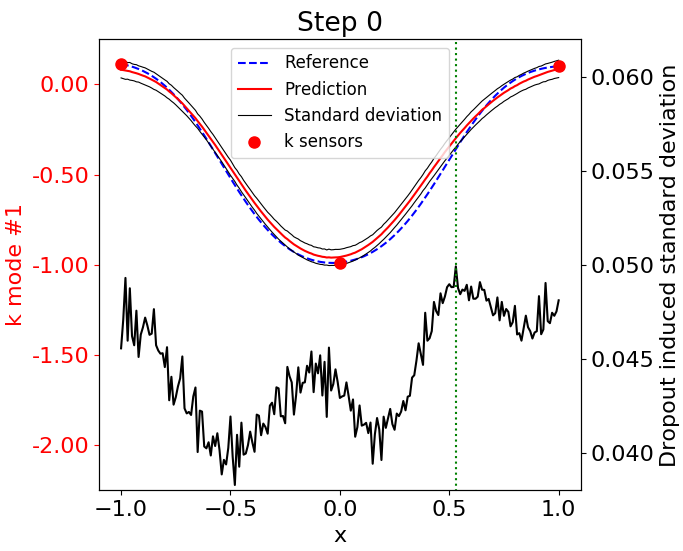}
		\caption{}
		\label{fig:43_kmode01}
	\end{subfigure}%
	\begin{subfigure}{0.33\textwidth}
		\centering
		\includegraphics[width=\linewidth]{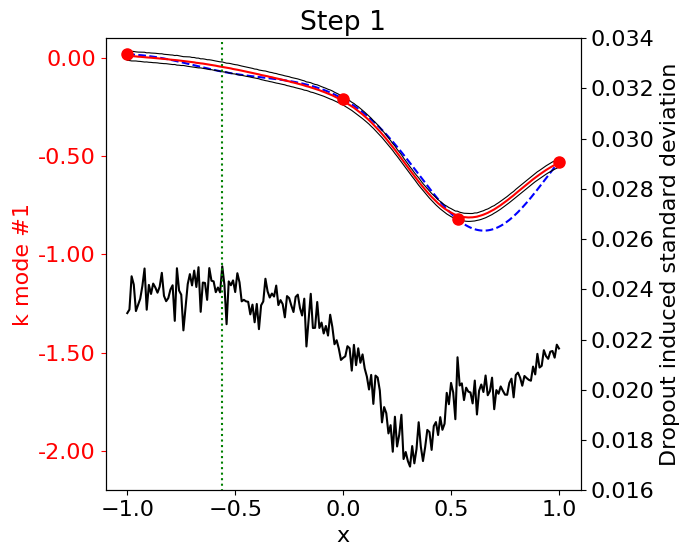}
		\caption{}
		\label{fig:43_kmode11}
	\end{subfigure}%
	\begin{subfigure}{0.33\textwidth}
		\centering
		\includegraphics[width=\linewidth]{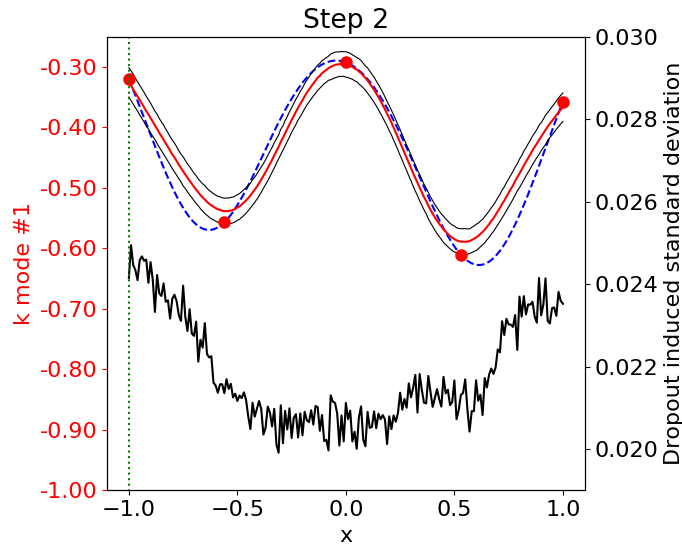}
		\caption{}
		\label{fig:43_kmode21}
	\end{subfigure}
	\caption{Active learning for stochastic inverse problems: The first modes of $k$ in the three steps and their associated standard deviation induced from the dropout neural network. The green dashed line indicates the location where the standard deviation reaches its maximum, and where the new sensor will be added in the next step.}
	\label{fig:43_kmodes1}
\end{figure}

\begin{figure}[htbp]
	\centering
	\begin{subfigure}{0.5\textwidth}
		\centering
		\includegraphics[width=\linewidth]{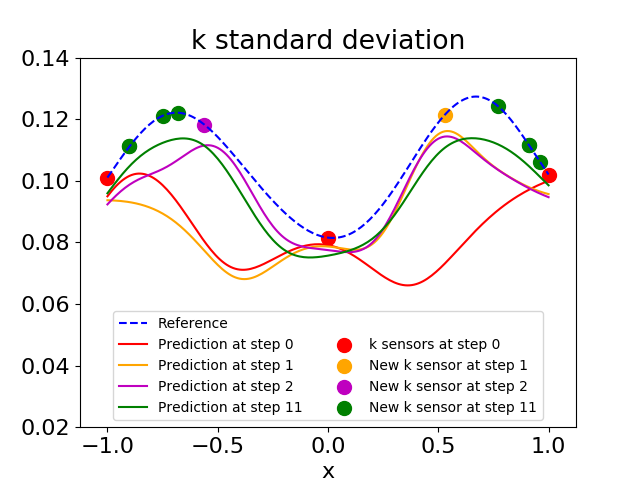}
		\caption{}
		\label{fig:43_kstd}
	\end{subfigure}%
	\begin{subfigure}{0.5\textwidth}
		\centering
		\includegraphics[width=\linewidth]{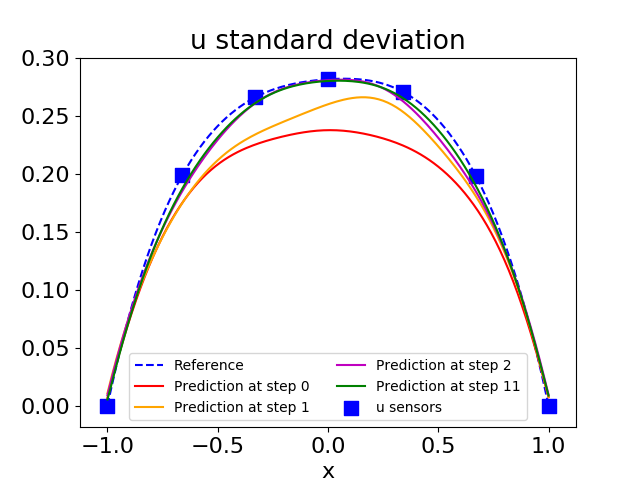}
		\caption{}
		\label{fig:43_ustd}
	\end{subfigure}
	\caption{Active learning for stochastic inverse problems: The predicted standard deviation of $k$ (in plot (a)) and $u$ (in plot (b)) in the first three and the last steps. In plot (a) the $k$-sensors at step 0 are colored in red and the newly added $k$-sensors are denoted in different colors. The $u$-sensors (as depicted in plot(b)) are kept the same. The reference solutions are calculated using the continuous trajectories that generate the training data.}
	\label{fig:43_std}
\end{figure}

We carry out the active learning steps until the $k$ prediction error does not decrease after a new sensor is added. To better illustrate the process of adding new sensors, in Figure \ref{fig:43_kmode01}, \ref{fig:43_kmode11} and \ref{fig:43_kmode21} we show the learned first mode of $k$ and its associated standard deviation from dropout uncertainty. Three steps of active learning are displayed. Note that the shapes of the first modes do not stay the same due to the fact that every time a new $k$-sensor is added, the principal components of $K$ in Equation \ref{eqn:pca} are therefore changing. Figure \ref{fig:43_kstd} and Figure \ref{fig:43_ustd} show the comparison of the predicted standard deviation of $k$ and $u$ in the first three steps and the last step, respectively. It is evident that adding extra $k$-sensors automatically according to the dropout-induced uncertainty will improve the accuracy of standard deviation prediction for both $k$ and $u$. Finally, the trained model is used to predict continuous trajectories of $k$ and $u$ in the test samples. We can conclude from Table \ref{tab:43_error} that adding extra $k$-sensors based on active learning with the dropout helps us to make better predictions.

\begin{table}[htbp]
\centering
\caption{Comparison of the relative $L_2$ error at different steps}
\label{tab:43_error}
\begin{tabular}{|c||c|c|c||c|c|c|}
\hline
        & $k$ mean & $k$ std   & $k$ prediction & $u$ mean & $u$ std   & $u$ prediction \\ \hline
Step 0  & $0.87\%$ & $26.57\%$ & $6.07\%$       & $0.18\%$ & $15.23\%$ & $3.06\%$       \\ \hline
Step 1  & $2.79\%$ & $20.19\%$ & $5.29\%$       & $0.24\%$ & $9.33\%$  & $2.41\%$       \\ \hline
Step 2  & $1.48\%$ & $10.21\%$ & $2.58\%$       & $0.14\%$ & $3.89\%$  & $1.08\%$       \\ \hline
Step 11 & $0.46\%$ & $8.49\%$  & $2.01\%$       & $0.04\%$ & $2.93\%$  & $0.67\%$       \\ \hline
\end{tabular}
\end{table}

\section{Summary}
\label{S:5}
We have presented a new approach to quantify the parametric uncertainty in the physics-informed neural networks (PINNs) by employing the arbitrary polynomial chaos (aPC) expansion to represent the stochastic solution. We use the data collected from sensor measurements to build a set of arbitrary polynomial basis and learn the modal functions of the aPC expansion through the PINNs, i.e., DNNs that encode the underlying stochastic differential equation. The proposed data-driven method can be used to solve forward problems, but more importantly, it deals with stochastic {\em inverse} problems. In the classical inverse problem, typically all the information available is for the solution, and we aim to identify the parameters. Here we selected to solve the inverse problems, where we have partial information available both for the solution and the parameter, which is a stochastic process. Once the model is trained with existing sensor data, i.e., historical data, it can be used to predict new instances of trajectories of the quantity of interest (solution or parameter) at very small additional computational cost.

We aim at quantifying two different types of uncertainties, i.e., the {\em parametric} uncertainty due to the stochastic equation, as well as the {\em approximation} uncertainty of the PINN. The latter represents how well the PINN is trained and how robust it is at the predicting stage. To this end, we adopt the dropout strategy for estimating the approximation uncertainty. Dropout is typically used to deal with over-fitting problems but it can also be exploited to quantify the approximation uncertainty at no additional cost. In our examples, we use DNNs to learn the modal functions of the stochastic parametric modes and the dropout strategy to quantify their associated uncertainty. Based on this, we propose an iterative method of actively learning where to place new sensors to enhance the approximation accuracy of PINNs. The numerical results exhibit the effectiveness of such an {\em active learning} strategy, not only in placing new sensors but also in making better use of the existing ones.

There are other possible methods of quantifying parametric and approximation uncertainties. For example, we can abandon aPC and use directly the stochastic data as the input. One possible approach is to consider the random space together with the physical space. Hence, standard DNNs, and in particular the deterministic PINNs developed in~\cite{MaziarParisGK17_1, MaziarParisGK17_2} can be directly used to stochastic differential equations. However, we did not choose this approach for two reasons: first, it does not lead to explicit expressions for stochasiticity of the quantity of interest (QoI); and second, different from sampling in the physical space where we can always mark the location by their coordinates, marking random instances with random variables is much harder. Moreover, the dimension of the random space is usually much higher than that of the physical space, requiring a proper dimension reduction procedure to be carried out before actually solving the problem. Although weaker, the high dimensionality is also a limitation of the aPC approach because it requires the high-dimensional DNN outputs, which could make the training process harder. Other promising methods that may be able to deal with high dimensional stochastic problems are the generative adversarial networks (GANs)~\cite{gans}, and we plan to systematically investigate this line of research in our future work.

\section{Acknowledgement}
This work is supported by DARPA N66001-15-2-4055, Air Force FA9550-17-1-0013, ARL W911NF-12-2-0023, NSF of China (No. 11671265) and the Science Challenge Project (No. TZ2018001). In addition, we would like to thank Liu Yang for his generous advice.

\section*{References}
\bibliography{bibfile}

\end{document}